\documentclass[12pt,reqno,twoside]{amsart}
\usepackage{amsmath,amsthm,amsfonts,amssymb,setspace,verbatim,cite}

\usepackage[english]{babel}
\usepackage[utf8]{inputenc}
\usepackage{times}
\usepackage[T1]{fontenc}
\usepackage{mathrsfs}
\usepackage{mathabx}
\usepackage{enumitem}

\usepackage[top=15mm,bottom=25mm,left=20mm,right=25mm,includeheadfoot,headsep=8mm,footnotesep=5mm]{geometry} 

\usepackage[hypertexnames=false]{hyperref}
\hypersetup{
pdfencoding=auto,
pageanchor=true,
plainpages=true,
pdftoolbar=true,
pdfmenubar=true,
pdffitwindow=false,
pdfstartview={FitH},
pdfnewwindow=true,
colorlinks=true,
linkcolor=black,
citecolor=black,
filecolor=black,
urlcolor=black,
linktoc=all,
breaklinks=true,
}

\newtheorem{theorem}{Theorem}[section]
\newtheorem{lemma}[theorem]{Lemma}
\newtheorem{corollary}[theorem]{Corollary}

\theoremstyle{definition}
\newtheorem{definition}[theorem]{Definition}
\newtheorem{hypothesis}[theorem]{Hypothesis}

\newtheorem{examplex}[theorem]{Example}
\newenvironment{example}
  {\pushQED{\qed}\examplex}
  {\popQED\endexamplex}

\theoremstyle{remark}
\newtheorem{remark}[theorem]{Remark}

\numberwithin{equation}{section}

\makeatletter
\def\th@plain{%
  \thm@notefont{\bfseries}% same as heading font
  \itshape % body font
}

\def\th@definition{%
  \normalfont % body font
  \thm@notefont{\bfseries}% same as heading font
}

\renewcommand{\thefootnote}{\fnsymbol{footnote}}
\newcommand{\CorAuth}
{Corresponding author.}

\newcommand{\makepapertitle}{
\pdfbookmark[1]{Title page}{Title_page}
\thispagestyle{empty}\vspace*{8mm}
\begin{spacing}{1.1}
\LARGE\MakeUppercase{\papertitle}
\end{spacing}\vspace*{5mm}
}

\newcommand{\paperfirstauthor}{%
{\large\sc\firstauthor}\\[2mm]%\footnote{\firstthanks}} \\[2mm]
{\rm\firstaddress \\
\firstemail}\\[8mm]
}

\newcommand{\papersecondauthor}{%
{\large\sc\secondauthor\footnote{\secondthanks}}\\[2mm]
{\rm\secondaddress \\
\secondemail}\\[8mm]
}

\newcommand{\makeabstract}{%
\begin{minipage}{0.9\textwidth}
{\small {\sc Abstract.}
 \paperabstract
}
\end{minipage}\vfill
}

\newcommand{\MakeFirstPageTwoAuthors}{
\begin{center}
  \makepapertitle
  \paperfirstauthor
  \papersecondauthor
  \makeabstract
  \begin{minipage}[t]{0.3\textwidth}
\raggedleft {\bf Date (revised final version):} 
  \end{minipage}\hspace{0.01\textwidth}
  \begin{minipage}[t]{0.6\textwidth}
    \today\\ (submitted on June 27, 2013; accepted on July 7, 2014)
  \end{minipage}\\[4mm]
\noindent%
  \begin{minipage}[t]{0.3\textwidth}
    \raggedleft {\bf Running head:} 
  \end{minipage}\hspace{0.01\textwidth}
  \begin{minipage}[t]{0.6\textwidth}
    \runninghead
  \end{minipage}\\[4mm]
\noindent%
  \begin{minipage}[t]{0.3\textwidth}
    \raggedleft {\bf How to cite:} 
  \end{minipage}\hspace{0.01\textwidth}
  \begin{minipage}[t]{0.6\textwidth}
    J. Math. Anal. Appl. {\bf 421} (2015), no.~1, 779--805.\\[1mm]
    \footnotesize{\url{http://dx.doi.org/10.1016/j.jmaa.2014.07.015}}
  \end{minipage}\\[4mm]
\noindent%
  \begin{minipage}[t]{0.3\textwidth}
    \raggedleft {\bf License:}  
  \end{minipage}\hspace{0.01\textwidth}
  \begin{minipage}[t]{0.63\textwidth} %\onehalfspacing  
    \copyright{} \the\year. This manuscript version is made available under the 
     \href{http://creativecommons.org/licenses/by-nc-nd/4.0/}{CC-BY-NC-ND 4.0 license}
  \end{minipage}\\
\bigskip    
\vfill
\end{center}

\renewcommand{\thefootnote}{}
\footnotetext[2]{2010 {\it Mathematics Subject Classification:} \thesubjclass}
\footnotetext[3]{{\it Key words and phrases:} \thekeywords}
\setcounter{footnote}{0}
\renewcommand{\thefootnote}{{\bf\,\alph{footnote}\alph{footnote}\alph{footnote}\,}}

\markboth{{\sc\firstauthor\ and \secondauthor}}{{\sc\runninghead}}

\setcounter{page}{0}
\newpage
}

\newcommand{\papertitle}%
{On discrete symplectic systems: Associated maximal and minimal linear relations and nonhomogeneous problems}

\newcommand{\runninghead}%
{On Discrete Symplectic Systems}

\newcommand{\firstauthor}%
{Stephen L. Clark}

\newcommand{\firstaddress}
{Department of Mathematics \& Statistics, 101 Rolla Building,\\
Missouri University of Science and Technology, Rolla, MO 65409-0020, USA}

\newcommand{\firstemail}%
{E-mail: sclark@mst.edu}

\newcommand{\secondauthor}%
{Petr Zem{\'{a}}nek}

\newcommand{\secondaddress}
{Department of Mathematics and Statistics, Faculty of Science, Masaryk University \\
Kotl{\'{a}}{\v{r}}sk{\'{a}} 2, CZ-61137 Brno, Czech Republic}

\newcommand{\secondemail}%
{E-mail: zemanekp@math.muni.cz}

\newcommand{\secondthanks}
{\CorAuth}

\newcommand{\paperabstract}%
{In this paper we characterize the definiteness of the discrete symplectic system, study a nonhomogeneous discrete
symplectic system, and introduce the minimal and maximal linear relations associated with these systems. Fundamental 
properties of the corresponding deficiency indices, including a relationship between the number of square summable 
solutions and the dimension of the defect subspace, are also derived. Moreover, a sufficient condition for the existence 
of a densely defined operator associated with the symplectic system is provided.
}

\newcommand{\thekeywords}%
{discrete symplectic system; time-reversed system; definiteness condition; nonhomogeneous problem; linear relation; 
deficiency index.}

\newcommand{\thesubjclass}%
{{\it Primary\/} 47A06; {\it Secondary\/} 39A70; 47B39; 39A12.}

\newcommand{\submittedto}%
{Journal of Mathematical Analysis and Applications}

\newcommand{\mmatrix}[1]{\left(\begin{matrix} #1
  \end{matrix}\right)}
\newcommand{\msmatrix}[1]{\left(\begin{smallmatrix} #1
  \end{smallmatrix}\right)}
\newcommand{\Ps}{\Psi}
\newcommand{\Psk}{\Ps_{k}}
\newcommand{\fk}{f_{k}}
\newcommand{\gk}{g_{k}}
\newcommand{\zk}{z_{k}}

\newcommand{\al}{\alpha}
\newcommand{\als}{\alpha^*}

\newcommand{\la}{\lambda}
\newcommand{\bla}{\bar{\lambda}}
\newcommand{\De}{\Delta}
\newcommand{\de}{\delta}
\newcommand{\Ga}{\Gamma}
\newcommand{\Ph}{\Phi}
\newcommand{\Phk}{\Phi_k}

\newcommand{\E}{\mathcal{E}}
\newcommand{\T}{\mathcal{T}}
\newcommand{\Tmax}{T_{\mathrm{max}}}

\newcommand{\Tmin}{T_{\mathrm{min}}}
\newcommand{\Tnod}{T_0}

\newcommand{\Zbb}{\mathbb{Z}}
\newcommand{\Rbb}{\mathbb{R}}
\newcommand{\Nbb}{\mathbb{N}}
\newcommand{\Cbb}{\mathbb{C}}
\newcommand{\Sbb}{\mathbb{S}}
\newcommand{\ellP}{\ell^{\,2}_{\Ps}(\oinftyZ)}
\newcommand{\ellPo}{{\ell}^{\,2}_{\Ps,1}(\oinftyZ)}
\newcommand{\ellPn}{\ell^{\,2}_{\Ps,0}(\oinftyZ)}

\newcommand{\tellPI}[1]{\tilde{\ell}^{\,2}_{\Ps}{#1}}
\newcommand{\tellP}{\tilde{\ell}^{\,2}_{\Ps}(\oinftyZ)}
\newcommand{\tellPo}{\tilde{\ell}^{\,2}_{\Ps,1}(\oinftyZ)}
\newcommand{\A}{\mathcal{A}}
\newcommand{\B}{\mathcal{B}}
\renewcommand{\C}{\mathcal{C}}
\newcommand{\D}{\mathcal{D}}
\newcommand{\K}{\mathcal{K}}
\newcommand{\V}{\mathcal{V}}
\newcommand{\J}{\mathcal{J}}
\newcommand{\mL}{\mathscr{L}}
\newcommand{\mI}{\mathcal{I}}

\newcommand{\mhI}{\widehat{\mI}}
\newcommand{\hilb}{\mathscr{H}}
\renewcommand{\S}{\mathcal{S}}

\newcommand{\X}{\mathcal{X}}
\newcommand{\Xchi}{\X}
\newcommand{\inner}[2]{\langle #1,\,#2 \rangle}
\newcommand{\innerP}[2]{\langle #1,\,#2 \rangle_{\Ps}}
\newcommand{\innerPI}[3]{\langle #2,\,#3 \rangle_{\Ps,#1}}
\newcommand{\normP}[1]{\|{#1}\|_{\Ps}}
\newcommand{\normPI}[2]{\|{#2}\|_{\Ps,#1}}
\newcommand{\norm}[1]{\|{#1}\|}
\newcommand{\abs}[1]{|{#1}|}

\newcommand{\oinftyZ}{\Nbb_0}%{[0,\infty)_{\sZbb}}
\newcommand{\sZbb}{{\scriptscriptstyle{\Zbb}}}
\DeclareMathOperator{\Ran}{ran}
\DeclareMathOperator{\Ker}{ker}
\DeclareMathOperator{\dom}{dom}
\DeclareMathOperator{\mul}{mul}
\DeclareMathOperator{\ran}{ran}
\DeclareMathOperator{\rank}{rank}

\DeclareMathOperator{\diag}{diag}
\DeclareMathOperator{\codim}{codim}
\DeclareMathOperator{\sgn}{sgn}

\newcommand{\W}{\mathcal{W}}
\newcommand{\qtextq}[1]{\quad\text{#1}\quad}
\newcommand{\qtext}[1]{\quad\text{#1 }\ }

\newcommand{\vp}{\varphi}

\newcommand{\dla}{d_{\la}}
\newcommand{\tdla}{\td_{\la}}
\newcommand{\td}{\tilde{d}}
\newcommand{\tz}{\tilde{z}}
\newcommand{\tw}{\tilde{w}}
\newcommand{\tu}{\tilde{u}}
\newcommand{\tv}{\tilde{v}}
\newcommand{\tZ}{\widetilde{Z}}
\newcommand{\tf}{\tilde{f}}
\newcommand{\tg}{\tilde{g}}
\newcommand{\hz}{\hat{z}}
\newcommand{\hy}{\hat{y}}
\newcommand{\ty}{\tilde{y}}
\newcommand{\tM}{\widetilde{M}}
\newcommand{\tMla}{\tM_\la}
\newcommand{\tMbla}{\tM_{\bla}}
\newcommand{\Mla}{M_\la}

\newcommand{\zkla}{z_{k}(\la)}
\newcommand{\zkkla}{z_{k+1}(\la)}

\newcommand{\zklas}{z_{k}^*(\la)}

\newcommand{\zjlas}{z_{j}^*(\la)}

\newcommand{\Sla}[1]{\text{\rm(S$_{#1}$})}
\newcommand{\Slaf}[2]{\text{\rm(S$_{#1}^{#2}$)}}

\hypersetup{
pdfencoding=unicode,
pdftitle={\papertitle},
pdfauthor={\firstauthor{} \& \secondauthor},
pdfsubject={Department of Mathematics \& Statistics, Missouri University of Science and Technology --- Department of
Mathematics and Statistics, Faculty of Science, Masaryk University},
pdfkeywords={\thekeywords},
pdfsubject=\paperabstract,
pdfencoding=auto
}

\begin{document}

%%%%%%%%%%%%%%%%%%%%%%%%%%%%%%%%%%%%%%%%%%%%%% FIRST PAGE %%%%%%%%%%%%%%%%%%%%%%%%%%%%%%%%%%%%%%%%%%%%%%%%%%%%%%%%%%%%

\MakeFirstPageTwoAuthors

%%%%%%%%%%%%%%%%%%%%%%%%%%%%%%%%%%%%%%%%%%% SECTION %%%%%%%%%%%%%%%%%%%%%%%%%%%%%%%%%%%%%%%%%%%%%%%%%%%%%%%%%%%%%%%%%%

\section{Introduction}\label{S:Intro}

The spectral theory for difference equations and systems, which include discrete analogs of Sturm--Liouville and 
Hamiltonian systems of differential equations, has a long history and a considerable literature which we will not attempt 
to delineate here other than to cite the following works, and the references therein, to give the reader a sense
of the scope of the subject in time and content, cf. 
\cite{slC.fG04,slC.fG.mZ07,slC.fG.mZ08,slC.fG.mZ09,gT00,bmS05:I,bmS05:II,jmB68,fvA64}. In connection with this, the 
development of a Weyl-Titchmarsh theory for discrete Hamiltonian and symplectic systems parallel to that which exists for 
Hamiltonian systems of differential equations is of relatively recent origin,  cf. 
\cite{mB.sS10,sjM.kmS07,slC.pZ10,rSH.pZ:W-TGLP,rSH.pZ:joint,gR.yS10,yS.hS11,hB13,hB.foN11:JDEA.a,hB.foN11:JDEA.b,slC.fG04}.
Our current paper contributes to this ongoing development.

We investigate the nonhomogeneous problem as well as the basic development of linear relations associated with discrete 
symplectic systems written in the so-called time-reversed form given by
 \begin{equation*}\label{Sla}\tag{S$_\la$}
   z_{k}(\la)=\Sbb_{k}(\la)\,z_{k+1}(\la) \quad\text{with }\ \Sbb_{k}(\la)\coloneq\S_{k}+\la \V_{k} \
                                               \text{ and }\ k\in\oinftyZ,
 \end{equation*}
where $\la\in\Cbb$ is the spectral parameter, $\oinftyZ\coloneq[0,\infty)\cap\Zbb$, and $\S_k$ and $\V_k$ are complex 
$2n\times 2n$ matrices such that
\begin{equation}\label{E:2}
  \S_k^*\J \S_k=\J,\quad \S_k^*\J\,\V_k\ \text{ is Hermitian}, \quad \V_k^*\J\,\V_k=0,  \qtextq{and}
  \Psk\coloneq\J\S_k\J\V_k^*\J\geq0,
\end{equation}
where $\J$ represents a $2n\times2n$ skew-symmetric matrix given by $\J\coloneq\msmatrix{\phantom{-}0 & I\\ -I & 0}$.

We note (viz. Lemma~\ref{l2.2}) that the first, second, and third identities in \eqref{E:2}, i.e., \eqref{E:2}(i)--(iii), 
can be combined into the single equality
 \begin{equation}\label{E:1.2}
   \Sbb_{k}^*(\bla)\,\J\,\Sbb_{k}(\la)=\J \quad\text{for all }\ \la\in\Cbb \ \text{and } k\in\oinftyZ,
 \end{equation}
where $\Sbb_{k}^*(\bla)\coloneq(\Sbb_{k}(\bla))^*$. Identity~\eqref{E:1.2} justifies the terminology \emph{symplectic
system} for \eqref{Sla}, though system~\eqref{Sla} corresponds to the well-known time-reversed discrete symplectic
system introduced in~\cite[Remark~4]{mB.oD97} only when $\la\in\Rbb$; particularly, the case when $\la=0$. In addition, 
system \eqref{Sla} can also be viewed as a perturbation of the original symplectic system $z_{k}=\S_k z_{k+1}$, i.e., of
\eqref{Sla} with $\la=0$, but for which the fundamental properties of symplectic systems remain true with appropriate, 
natural, modifications. 

In~\cite{mB.sS10,slC.pZ10}, the Weyl--Titchmarsh theory was first established for discrete symplectic systems given by 
 \begin{equation}\label{E:1.1.1}
   z_{k+1}(\la)=\Sbb_{k}^{-1}\,z_k(\la), \quad k\in\Nbb_0,
 \end{equation}
in which a special form for $\V_k$ is assumed; the proper generalization is later derived in \cite{rSH.pZ:W-TGLP}, 
see also \cite{rSH.pZ:ICDEA}. The results given in \cite{rSH.pZ:W-TGLP} for system~\eqref{E:1.1.1} remain 
valid for system~\eqref{Sla} with standard changes given for the definition of the semi-inner product, 
viz. \eqref{E:inner.prod} (cf. \cite[Theorem~2.8 and Section~4]{rSH.pZ:W-TGLP}), and for the associated weight 
function, viz. \eqref{E:2}(iv) (cf.\cite[Identity~(1.1)(iv)]{rSH.pZ:W-TGLP}). 

Consideration here of the time-reversed form given in system~\eqref{Sla}, rather than that given in system~\eqref{E:1.1.1}, 
is motivated, in part, by a~desire to produce more natural calculations involving the semi-inner product and in particular 
a~more natural form for a~Green function associated with nonhomogeneous discrete symplectic systems, 
viz. Lemma~\ref{L:green.prop}. We can also associate with system~\eqref{Sla} a~densely defined operator, because there is 
no shift in the associated semi-inner product (cf. Theorem~\ref{T:dense.domain} and \cite{gR14:AML}). Moreover, this 
approach will enable us in subsequent research to generalize these results and unify them with the continuous time case 
by means of the time scale theory.

Given the inherent semi-definiteness of the function $\Psi$  defined in \eqref{E:2} (cf. \eqref{e2.2}), it is natural to
consider the construction of \emph{linear relations} in association with \eqref{Sla}, their extensions and their associated
spectral theory. The theory of linear relations provides powerful tools for the study of \emph{multivalued} linear 
operators in a Hilbert space, especially for non-densely defined linear operators. The study of linear relations in this
context traces back to \cite{rA61}; see also \cite{sH.hsvdS.fhS09,eaC73:Ex,rC98,aD.hsvdS74} and 
the references therein. For linear Hamiltonian differential systems given by
 \begin{equation}\label{E:cLHS}
   -\J z'(t)=\big[H(t)+\la W(t)\big] z(t),
 \end{equation}
where $H(t)$ and $W(t)$ are Hermitian and $W(t)$ is positive semi-definite, this approach was initiated in~\cite{bcO69} and 
further developed, e.g., in \cite{mL.mmM03,jB.sH.hsvdS.rW11,isK83,sH.hsvdS.hW00}. 

In \cite{gR.yS11}, linear relations are considered in association with the linear Hamiltonian difference system 
 \begin{equation}\label{E:dLHS}
   \De\!\mmatrix{x_k\\ u_k}=\big(H_k+\la W_k\big)\mmatrix{x_{k+1}\\ u_k}, \quad
   H_k\coloneq\mmatrix{A_k & B_k\\ C_k & -A_k^*}, \quad W_k\coloneq\mmatrix{\phantom{-}0 & W^{[1]}_k\\ -W^{[2]}_k & 0},
 \end{equation}
where $B_k$ and $C_k$ are Hermitian matrices, $W^{[j]}_k\geq 0$, $j=1,2$, and the matrix $\tilde{A}_k\coloneq I-A_k$ is
invertible. Here, we note that the invertibility assumption and the additional requirement 
$W^{[1]}_k(I-A_k)^{-1}W^{[2]}_k\equiv0$ imply that system~\eqref{E:dLHS} can be written as a discrete system whose form is 
given by~\eqref{Sla}; cf. \cite[Formula~(2.3)]{yS06}. On the other hand, with the supplementary assumption concerning the
invertibility of the $n\times n$ matrix in the left upper block of the matrix $\Sbb_k(\la)$, $k\in\oinftyZ$, system
\eqref{Sla} can be written as the linear Hamiltonian difference system but with a \emph{nonlinear} dependence on the 
spectral parameter, see \cite{rSH.pZ:ICDEA}. Moreover, using the time scale theory, see e.g. \cite{oD.rH01}, it can be 
seen that the discrete symplectic systems given by \eqref{Sla} provide a proper discrete analogue of linear Hamiltonian 
differential systems. Finally, let us note that system~{\eqref{Sla}} includes any even order Sturm--Liouville difference 
equation; cf. \cite[Remark~2]{mB.oD97}, \cite{pZ13}, and see also Example~\ref{Ex:definit}.

Hence, we shall introduce minimal and maximal linear relations associated with our discrete symplectic system and establish 
fundamental properties for them in analogy with \cite[Section~2]{mL.mmM03} for system~\eqref{E:cLHS} and 
\cite[Section~5]{gR.yS11} for system~\eqref{E:dLHS}. Moreover, the reader can observe an essential difference, one which 
appears natural in the context of the time scale theory, in the assumptions for systems \eqref{E:cLHS} and \eqref{Sla} 
concerning the invertibility of $W$ and $\Ps$, respectively. While the matrix $W$ can be invertible, the matrix $\Psk$ is 
singular for every $k\in\oinftyZ$, see \eqref{E:2}.

The remainder of this paper is organized as follows. In Section~\ref{S:2}, we present fundamental properties of
system~\eqref{Sla} and recall some basic facts from the theory of linear relations. In Section~\ref{S:3} we discuss the
definiteness condition for system~\eqref{Sla}, which plays a crucial role in the spectral theory, and derive some equivalent
characterizations. A nonhomogeneous discrete symplectic system is studied in Section~\ref{S:nonhom}. Concluding with
Section~\ref{S:4}, the maximal and minimal linear relations associated with the discrete symplectic systems are introduced
and their fundamental properties, such as a relationship between the deficiency indices of the minimal relation in a suitable
Hilbert space and the number of square summable solutions of system \eqref{Sla}, are established. In this final section, we
also present a~sufficient condition providing the existence of a densely defined operator associated with
the discrete symplectic system. 

%%%%%%%%%%%%%%%%%%%%%%%%%%%%%%%%%%%%%%%%%%% SECTION %%%%%%%%%%%%%%%%%%%%%%%%%%%%%%%%%%%%%%%%%%%%%%%%%%%%%%%%%%%%%%%%%%

\section{Preliminaries}\label{S:2}
\subsection{Notation}\label{S:2.0}

Throughout this paper, matrices will be considered over the field of complex numbers $\Cbb$. For any $\la\in\Cbb$ we 
denote its imaginary part by $\Im(\la)$. By $\Cbb^{r\times s}$, $r,s\in\Nbb$, we mean the space of $r\times s$ complex 
matrices and $\Cbb^{r\times 1}$ will be denoted by $\Cbb^r$ for $r\in\Nbb$. With $M\in\Cbb^{r\times s}$, let $M^\top$ denote 
the transpose, and let $M^*$ denote the adjoint or conjugate transpose of the matrix $M$; for parameter dependent matrices, 
$M^*(\lambda)\coloneq M(\lambda)^*$. Let $M\ge 0$ and $M\le 0$ indicate that $M$ is positive or negative semi-definite, 
respectively. Similarly, $M>0$ (respectively, $M<0$) denotes a positive definite (respectively, negative definite) matrix. 
By $\J$ we denote the real $2n\times 2n$ skew-symmetric matrix given as
 \begin{equation}\label{e2.1}
  \J \coloneq\left(\hspace{-4pt}\begin{array}{cc} 0& I_n\\-I_n&0 \end{array}\!\right),
  \end{equation}
where $I_n$ is the $n\times n$ identity matrix.

If $\mI$ denotes an interval in $\Rbb$, then the associated discrete interval in the set of integers, $\Zbb$, is denoted 
by $\mI_\Zbb\coloneq \mI\cap\Zbb$. In particular, $\Nbb=[1,\infty)_\Zbb$, and $\oinftyZ=[0,\infty)_\Zbb$. However, in
practice, $\mI_\Zbb$ will be referred as the \emph{discrete interval} $\mI$. Hence, for the discrete interval $\mI$, 
by $\Cbb(\mI)^{r\times s}$ we denote the space of sequences, defined on $\mI$, of complex $r\times s$ matrices, where 
typically $r\in\{n, 2n\}$ and $1\le s\le 2n$. If $\S\in\Cbb(\mI)^{r\times s}$, then $\S(k)\coloneq\S_k$ for $k\in \mI$;
if $\S(\lambda)\in\Cbb(\mI)^{r\times s}$, then $\S(\lambda,k)\coloneq\S_k(\lambda)$ for $k\in \mI$. When 
$\S\in\Cbb(\mI)^{m\times s}$, and $\T\in\Cbb(\mI)^{s\times n}$, then $\S\T\in\Cbb(\mI)^{m\times n}$, where 
$(\S\T)_k\coloneq\S_k\T_k,\ k\in\mI$. The set $\Cbb_0(\mI)^{r\times s}$ represents the subspace of $\Cbb(\mI)^{r\times s}$ 
consisting of sequences compactly supported in the discrete interval $\mI$. The symbol $\Delta$ means the forward difference 
operator acting on $\Cbb(\mI)^{r\times s}$, where $(\Delta z)_k\coloneq z_{k+1} - z_k$, for all $k\in \mI$ and all 
$z\in\Cbb(\mI)^{r\times s}$. We shall also use the customary equivalence given by $(\Delta z)_k\coloneq \Delta z_k$.  
Moreover, we let $z_k \big|_{m}^n\coloneq z_n - z_m$.

\subsection{Time-reversed discrete symplectic systems}\label{S:2.1}

The basic relations given in the next lemma are easily shown and used throughout.

\begin{lemma}\label{l2.1}
 The following is true with $\J$ defined as in \eqref{e2.1}, and $\S,\V,\Psi\in\Cbb^{2n\times 2n}:$
  \begin{itemize}
   \item[\rm{(i)}] $\S^*\J\S=\J$ if and only if $\S\J \S^*=\J$;
   \item[\rm{(ii)}] $\S^*\J\S=\J$ if and only if $\S^{-1}=-\J\S^*\J$;
   \item[\rm{(iii)}] if $\V\coloneq\J^*\Psi\S=-\J\Psi\S$, where $\S^*\J\S=\J$, and
                      \begin{equation}\label{e2.2}
                       \Psi^*=\Psi, \quad \Psi\J\Psi=0,
                      \end{equation}
                     then
                      \begin{equation}\label{e2.3}
                       (\V^*\J\S)^*=\V^*\J\S, \quad (\V\J\S^*)^*=\V\J\S^*,\quad \V^*\J\V=\V\J\V^*=0,
                      \end{equation}
                     and $\Psi=\J\S\J\V^*\J$.  However, if $\S^*\J\S =\J$ and $\V\in\Cbb^{2n\times 2n}$ satisfies
                     \eqref{e2.3}, then $\Psi\coloneq\J\S\J\V^*\J$ satisfies  \eqref{e2.2}, and $\V=\J^*\Psi\S$.
  \end{itemize}
\end{lemma}

Lemma~\ref{l2.1} establishes a correspondence between the matrix pairs $\{ \S, \Psi\}$ and $\{\S,\V\}$ in which $\S$
satisfies $\S^*\J\S=\J$, $\Psi$ satisfies \eqref{e2.2} and $\V$ satisfies \eqref{e2.3}. The first part of the next result is
also easily verified, while the second part follows as in the proof of \cite[Lemma~2.2]{rSH.pZ:W-TGLP}. Let us note a typo 
in the latter reference, where the absolute value is missing.

\begin{lemma}\label{l2.2}
 Let $\Sbb(\lambda)= \S + \lambda\V$, where $\S,\V\in\Cbb^{2n\times 2n}$, $\lambda\in\Cbb$. Then
 $\Sbb^*(\bar\lambda)\J\Sbb(\lambda) = \J$  for all $\lambda\in\Cbb$ if and only if $\S^*\J\S=\J$, and $\V$ satisfies 
 \eqref{e2.3}. Moreover, $\abs{\det \Sbb(\la)}=1$.
\end{lemma}

Relevant to these two lemmas, we assume that the next conditions hold for the remainder of the paper.

\begin{hypothesis}\label{h2.3}
 For $\S, \Psi\in\Cbb(\oinftyZ)^{2n\times 2n}$, and for each $k\in\oinftyZ$, the following is assumed:
  \begin{enumerate}
   \item $\S_k^*\J\S_k=\J$;
   \item $\Psi_k$ satisfies \eqref{e2.2} and  $\Psi_k\ge 0$;
   \item $\Sbb_k(\lambda) \coloneq \S_k + \lambda\V_k$, where $\V_k \coloneq -\J\Psi_k\S_k$.
  \end{enumerate}
\end{hypothesis}
Given this hypothesis, note that $\V_k$ satisfies \eqref{e2.3} for all  $k\in\oinftyZ$, then by Lemma~\ref{l2.1} 
we have
 \begin{equation*}
  \Psi_k= \J\S_k\J\V^*_k\J,\quad k\in\oinftyZ,
 \end{equation*}
and by Lemma~\ref{l2.2} also $\Sbb_k^*(\bar\lambda)\J\Sbb_k(\lambda)=\J$ for all $k\in\oinftyZ$ and all 
$\lambda\in\Cbb$.

The discrete symplectic system in time-reversed form is now given by
 \begin{equation*}\tag{S$_\la$}
   z_{k}(\la)=\Sbb_{k}(\la)\,z_{k+1}(\la), \quad k\in\oinftyZ,
 \end{equation*}
where $z(\lambda)\in\Cbb(\oinftyZ)^{2n\times m}$, and $1\le m \le 2n$. In the future, \Sla{\nu} will denote a system of 
the form given in \eqref{Sla} with the parameter $\la$ replaced by $\nu$. Note that identity~{\eqref{E:1.2}}, the 
invertibility of {$\Sbb_k(\lambda)$} from Lemma~{\ref{l2.2}}, and $\J^{-1}=-\J$ imply, for all $k\in\oinftyZ$ and all 
$\lambda\in\Cbb$, that
 \begin{equation}\label{E:1A}
  \Sbb_k^{-1}(\la)=-\J\,\Sbb^*_k(\bla)\,\J,
 \end{equation}
with the consequent existence of unique solutions on $\oinftyZ$, $z(\lambda)\in\Cbb^{2n\times m}(\oinftyZ)$, for system
$\Sla{\lambda}$.

\begin{lemma}[Wronskian-type identity]\label{L:3}
 Let $\la\in\Cbb$ and $m\in\Nbb$. For solutions $z(\la)\in\Cbb(\oinftyZ)^{2n\times m}$ and 
 $z(\bla)\in\Cbb(\oinftyZ)^{2n\times m}$ of \eqref{Sla} and \Sla{\bla}, respectively, we have, for all $k\in\oinftyZ$,
  \begin{equation*}\label{E:9}
    z_{k}^*(\la)\,\J z_{k}(\bla)=z_{0}^*(\la)\,\J z_{0}(\bla).
  \end{equation*}
\end{lemma}
\begin{proof}
 This follows directly from Lemma~\ref{l2.2} since
  \begin{equation*}
   z_{k}^*(\la)\,\J z_{k}(\bla)
    =z_{k+1}^*(\la)\,\Sbb_k^*(\la)\,\J\,\Sbb_k(\bla)\,z_{k+1}(\bla)
    =z_{k+1}^*(\la)\,\J z_{k+1}(\bla),
  \end{equation*}
 for all $k\in\oinftyZ$, $\la\in\Cbb$.
\end{proof}

Throughout, we shall let $\Ph(\la)\in\Cbb(\oinftyZ)^{2n\times 2n}$ denote a fundamental system of solutions for~\eqref{Sla}.  
If this fundamental system is such that, for some $k_0\in\oinftyZ$,
  \begin{equation}\label{E:2.4}
    \Ph^*_{k_0}(\la)\,\J\,\Ph_{k_0}(\bla)=\J,
  \end{equation}
then, as an immediate consequence of the preceding lemmas,
  \begin{equation}\label{E:18A}
    \Phk^*(\la)\,\J\,\Phk(\bla)=\J,\quad
    \Ph^{-1}_k(\la)=-\J\,\Phk^{*}(\bla)\,\J,\quad
    \Phk(\la)\,\J\,\Phk^*(\bla)=\J,
  \end{equation}
for all $k\in\oinftyZ$.  When $\Phi_{k_0}(\la)=\Phi_{k_0}(\bla)=C\in\Cbb^{2n\times 2n}$, note that \eqref{E:2.4} simply 
states that $C$ is a~symplectic matrix.

We associate with the homogeneous system~\eqref{Sla} the following nonhomogeneous system
 \begin{equation}\label{Slaf}\tag{S$_\la^f$}
   \zkla=\Sbb_k(\la) \zkkla-\J\Psk\fk,\quad k\in\Nbb_0,
\end{equation}
where $f\in\Cbb(\oinftyZ)^{2n\times m}$, $1\le m\le 2n$. In analogy to previous notation, \Slaf{\nu}{g} denotes the 
nonhomogeneous system of the form given in \eqref{Slaf}, but with $\la$ replaced by $\nu$ and $f$ replaced by $g$. When 
convenient, we shall suppress the dependence of $z$ on $\la$ when $\la=0$. Note also that \eqref{Sla} is equivalent to 
\Slaf{\la}{0}.

The next result presents an essential tool used throughout; cf. \cite[Theorem~2.6]{rSH.pZ:W-TGLP} and 
\cite[Identity~(5.3)]{slC.pZ10}.

\begin{theorem}[Extended Lagrange identity]\label{T:1}
 Let $\la,\nu\in\Cbb$ and $1\le m\le 2n$. If $z(\la)\in\Cbb(\oinftyZ)^{2n\times m}$ and 
 $u(\nu)\in\Cbb(\oinftyZ)^{2n\times m}$ are solutions of systems~\eqref{Slaf} and \Slaf{\nu}{g}, respectively, with 
 $f,g\in\Cbb(\oinftyZ)^{2n\times m}$.  Then,
  \begin{gather}\label{E:7}
   \De\big[\zklas\,\J u_{k}(\nu) \big]=(\bar{\la}-\nu)\,\zklas\,\Psk\,u_{k}(\nu)
                                 +f^*_k\,\Psk\,u_{k}(\nu)-\zklas\,\Psk\,\gk,\\
   \zjlas\,\J u_j(\nu) \Big|_{0}^{k+1}=(\bar{\la}-\nu)\sum_{j=0}^{k}\zjlas\,\Ps_j\, u_j(\nu)
                                       +\sum_{j=0}^{k}f^*_j\,\Ps_j\, u_j(\nu)-\sum_{j=0}^{k}\zjlas\,\Ps_j\,g_j. \label{E:8}
  \end{gather}
\end{theorem}
\begin{proof}
By Lemma~\ref{l2.1} (cf. \eqref{E:2}) and \eqref{E:1A}, we see that
  \begin{align*}
   \De [&\zklas\,\J u_{k}(\nu) ]\\
    &=\big[\Sbb^{-1}_k(\la)\,\zkla+\Sbb^{-1}_k(\la)\,\J\,\Psk\fk\big]^*\J
      \big[\Sbb^{-1}_k(\nu)\,u_{k}(\nu) +\Sbb^{-1}_k(\nu)\,\J\,\Psk\gk\big]-\zklas\,\J\,u_{k}(\nu) \\
    &=\zklas\,\big[\Sbb_k^{*-1}(\la)\,\J\,\Sbb_k^{-1}(\nu)-\J\big]\, u_{k}(\nu)
      +f^*_k\,\Psk\,u_{k}(\nu)-\zklas\,\Psk\,\gk \\
    &=(\bar{\la}-\nu)\,\zklas\,\Psk\,u_{k}(\nu) +f^*_k\,\Psk\, u_k(\nu)-\zklas\,\Psk\,\gk,
  \end{align*}
 which yields \eqref{E:7}; which in turn, by summation, yields \eqref{E:8}.
\end{proof}

Finally, the following equivalent expression of system \eqref{Slaf} will play a key role in the results established in
Section~\ref{S:4}.

\begin{lemma} \label{L:slaf.equiv}
 System \eqref{Slaf} can be written equivalently as
  \begin{equation} \label{E:nonhom.equiv}
    \J\big(\zkla-\S_k\zkkla\big)=\la\Psk\zkla+\Psk\fk,\quad k\in\Nbb_0.
  \end{equation}
\end{lemma}
\begin{proof}
 Since one can easily observe from the identities in~\eqref{e2.3} that it holds
  \begin{equation}\label{E:slaf.equiv.proof.1}
   \Sbb_k(\la)=(I-\la\,\J\,\Psk)\,\S_k \qtextq{for all }\ k\in\Nbb_0,
  \end{equation}
 system \eqref{Slaf} can be also expressed as 
  \begin{equation}\label{E:slaf.equiv.proof.2}
   z_k(\la)=(I-\la\,\J\,\Psk)\,\S_k\,z_{k+1}(\la)-\J\,\Psk\,f_k,\quad k\in\Nbb_0.
  \end{equation}
 Since $(I-\la\,\J\,\Psk)^{-1}=(I+\la\,\J\,\Psk)$, which follows immediately from {\eqref{e2.2}}, we obtain from 
 \eqref{E:slaf.equiv.proof.2} that
  \begin{equation*}
   \S_k z_{k+1}(\la)=z_k(\la)+\la\,\J\,\Psk\,z_k(\la)+\J\,\Psk\,f_k, \quad k\in\Nbb_0,
  \end{equation*}
which after multiplication by the matrix $\J$ from the left yields \eqref{E:nonhom.equiv}.
\end{proof}

If we denote the left side of \eqref{E:nonhom.equiv} by
 \begin{subequations}\label{e2.9}
  \begin{equation}
   \mL(z(\la))_k\coloneq\J(\zkla-\S_k\zkkla),
  \end{equation}
then it follows from Lemma~\ref{L:slaf.equiv} that system \eqref{Slaf} is equivalent to
 \begin{equation}
   \mL(z(\la))_k=\la\Psk\zkla+\Psk\fk,\quad k\in\Nbb_0.
 \end{equation}
\end{subequations}
Hence $\mL$ represents a  linear map on $\Cbb(\oinftyZ)^{2n\times m}$  with  $1\le m\le 2n$, and \eqref{e2.9} will be 
summarized by writing $\mL(z(\la)) = \la\Psi z(\la) + \Psi f$.

%%%%%%%%%%%%%%%%%%%%%%%%%%%%%%%%%%%%%%%%%%%%%%%%%%%%%
\subsection{Sequence spaces}\label{SS:2.2}

With respect to the assumption in Hypothesis~\ref{h2.3} that $\Psi\in\Cbb(\oinftyZ)^{2n\times 2n}$
is a~sequence with positive semi-definite terms, we define a semi-inner product for $\Cbb(\oinftyZ)^{2n}$ by
 \begin{equation}\label{E:inner.prod}
  \innerPI{\mI}{z}{w}\coloneq\sum_{k\in\mI}\zk^*\,\Psk\,w_k,
 \end{equation}
where $z, w\in\Cbb(\oinftyZ)^{2n}$, and $\mI\subseteq\Nbb_0$ is a discrete subinterval. The associated linear space 
of all square summable sequences is denoted by
 \begin{equation*}
  \ell^2_\Psi(\mI)\coloneq \big\{z\in\Cbb(\oinftyZ)^{2n}\mid\ \normPI{\mI}{z}<\infty\big\},
 \end{equation*}
where $\normPI{\mI}{\cdot}\coloneq \sqrt{\innerPI{\mI}{\cdot}{\cdot}}$ denotes the associated semi-norm. When 
$\mI=\Nbb_0$ the corresponding semi-inner product and semi-norm will be denoted by $\innerP{\cdot}{\cdot}$, and 
$\normP{\cdot}$, respectively.
The quotient space obtained by factoring out the kernel of the semi-norm is denoted by
 \begin{equation*}
  \tilde\ell^2_\Psi(\mI)\coloneq\ell^2_\Psi(\mI)\big/\big\{z\in\Cbb(\Nbb_0)^{2n}\mid \
  \normPI{\mI}{z}=0\big\}.
 \end{equation*}
Then, with respect to the equivalence classes $\tz\in\tilde\ell^2_\Psi(\mI)$, we see that $\tilde\ell^2_\Psi(\mI)$ 
is a Banach space with respect to the norm generated by the quotient space map, $\pi(z)=\tz$, and a Hilbert space with
respect to the associated inner product where $\innerPI{\mI}{\tz}{\tw} \coloneq \innerPI{\mI}{z}{w} $ 
(cf. \cite[Lemma~2.5]{yS06}).

In light of this notation, we define
 \begin{equation*}\label{E:4.3}
   \ellPn\coloneq \big\{z\in \Cbb_0(\Nbb_0)^{2n}\mid\  z_0=0 \big\},
 \end{equation*}
and 
 \begin{equation*}
  \ellPo\coloneq \big\{z\in\Cbb(\Nbb_0)^{2n}\cap\ell^2_\Psi(\mI)\mid\ \exists\ N\in\oinftyZ\ \text{such that} \
                                                 \Psk\,\zk=0 \text{ for }k\geq N\big\}.
 \end{equation*}

%%%%%%%%%%%%%%%%%%%%%%%%%%%%%%%%%%%%%%%%%%%%%%%%
\subsection{Linear relations}\label{S:2.2}

We now recall some basic facts from the theory of linear relations relevant to the case at hand; cf.
\cite{sH.hsvdS.fhS09,eaC73:Ex,rC98,aD.hsvdS74}. Let $\hilb$ denote a Hilbert space over the field of complex numbers, 
$\Cbb$, with an inner product $\inner{\cdot}{\cdot}$. A (\emph{closed}) \emph{linear relation} $T$ in $\hilb$ is a (closed)
linear subspace of the product space $\hilb^2\coloneq \hilb\times\hilb$, i.e., the Hilbert space of all ordered  pairs
$\{z,f\}$ such that $z,f\in\hilb$. The \emph{domain}, \emph{range}, \emph{kernel}, and the \emph{multivalued part} of $T$ 
are respectively defined as
 \begin{gather*}
   \dom T\coloneq \big\{z\in\hilb\mid \exists f\in\hilb, \{z,f\}\in T\big\},\quad
   \ran T\coloneq \big\{f\in\hilb\mid \exists z\in\hilb,  \{z,f\}\in T\big\},\\
   \ker T\coloneq \big\{z\in\hilb\mid \{z,0\}\in T\big\},\quad
   \mul T\coloneq \big\{f\in\hilb\mid \{0,f\}\in T\big\}.
 \end{gather*}

In general, we let $T(z)\coloneq\{f\in\hilb\mid \{z,f\}\in T\}$, and note that a linear relation, $T$, is the graph of a 
linear operator in $\hilb$ when $T(0)=\{0\}$,  i.e., when the subspace $\mul T$ is trivial. The \emph{inverse} of $T$,
denoted as $T^{-1}$, is the linear relation
 \begin{equation*}
  T^{-1}\coloneq \big\{\{f,z\}\mid\ \{z,f\}\in T\big\}
 \end{equation*}
and note that $\dom T^{-1}=\ran T$, $\ran T^{-1}=\dom T$, $\ker T^{-1}=\mul T$, and $\mul T^{-1}=\ker T$. The 
\emph{adjoint} $T^*$ of the linear relation $T$ is defined by
 \begin{equation}\label{E:2.2.1A}
  T^*\coloneq \big\{\{y,g\}\in \hilb^2\mid\ \inner{z}{g}=\inner{f}{y},\ \forall\ \{z,f\}\in T\big\}.
 \end{equation}
The definition of $T^*$ reduces to the standard definition for the graph of the adjoint operator when $T$ is a~densely 
defined operator. Then $T^*$ is a closed linear relation and we have
 \begin{equation}\label{E:2.2.1}
  T^*=(\overline{T})^*,\quad
  T^{**}=\overline{T},\quad
  \ker T^*=(\ran T)^{\bot}=(\ran \overline{T})^\bot,\qtextq{and}
  (\dom T)^\bot=\mul T^*,
 \end{equation}
where $\overline{T}$ denotes the closure of $T$.  A linear relation $T$ is said to be \emph{symmetric} (or Hermitian) if 
$T\subseteq T^*$, and it is said to be \emph{self-adjoint} if $T^*=T$. It is easily seen that $T$ is a symmetric linear 
relation if and only if $\inner{z}{g}=\inner{f}{y}$ for all $\{z,f\},\{y,g\}\in T$.

With $\la\in\Cbb$, and with $T$ a linear relation, we define the linear relation, $T-\la I$, by
 \begin{equation}\label{e2.18a}
  T-\la I\coloneq \big\{\{z,f-\la z\}\mid\ \{z,f\}\in T\big\},
 \end{equation}
and note that $(T-\la I)^*=T^*-\bar\la I$. Then,  
 \begin{equation}\label{e2.17}
  \Mla(T)\coloneq \ran(T-\bla I)^\bot
 \end{equation}
is said to be the \emph{defect subspace} of $T$ and $\la$, and its dimension, i.e., 
 \begin{equation}\label{e2.18}
  d_{\la}(T)\coloneq \dim(\ran(T-\bla I)^\bot),
 \end{equation}
is said to be the \emph{deficiency index} of $T$ and $\la$. Since
 \begin{equation*}\label{E:2.2.2}
  \ran(T-\bla I)^\bot=\ker(T^*-\la I)=\big\{z\in \hilb\mid\ \{z,\la z\}\in T^*\big\},
 \end{equation*}
the deficiency indices of $T$ and $\overline{T}$ with the same $\la$ are equal by~\eqref{E:2.2.1}; 
cf. \cite[Lemma~2.4]{yS12}. We let
 \begin{equation*}
  d_+(T)\coloneq d_{i}(T),\qquad d_-(T)\coloneq d_{-i}(T),
 \end{equation*}
denote the so-called \emph{positive} and \emph{negative} deficiency indices of $T$, respectively. If $T$ is a symmetric 
linear relation, the values of $\dla(T)$ are constant in the open upper and lower half-planes of $\Cbb$; cf. 
\cite[Theorem~2.13]{yS12}. The linear relation $T$ has self-adjoint extensions if and only if the positive and negative 
deficiency indices are equal; cf. \cite[Corollary~6.4]{aD.hsvdS74}. Finally, for a closed symmetric linear relation $T$, 
it was shown in \cite[Lemma~2.25]{mL.mmM03}, when $\la\in\Rbb$ and $\ker(T-\la I)=\{0\}$, that
 \begin{equation}\label{E:2.2.3}
  d_\la(T)\leq d_\pm(T).
 \end{equation}

%%%%%%%%%%%%%%%%%%%%%%%%%%%%%%%%%%%%%%%%%%% SECTION %%%%%%%%%%%%%%%%%%%%%%%%%%%%%%%%%%%%%%%%%%%%%%%%%%%%%%%%%%%%%%%%%%

\section{Definiteness}\label{S:3}

In this section, we characterize the definiteness condition associated with the semi-inner product
\eqref{E:inner.prod} when the sequence, $\Psi\in\Cbb(\oinftyZ)^{2n\times 2n}$, possesses positive semi-definite
elements. This condition plays a significant role in the spectral theory, particularly the Weyl--Titchmarsh theory,
associated with the linear expression $\mL$ given in \eqref{e2.9}; cf. \cite{rSH.pZ:W-TGLP}. The condition is 
frequently called as the \emph{Atkinson condition}; cf.~\cite[Inequality~(3.7.10)]{fvA64}. Similar treatment in connection 
with the linear Hamiltonian differential and difference systems can be found in \cite{mL.mmM03,jB.sH.hsvdS.rW11} and \cite{gR.yS11}, 
respectively.

\begin{definition}\label{D:1}
 System~\eqref{Sla} is said to be \emph{definite} on a nonempty discrete interval $\mI\subseteq\oinftyZ$ if,
 for each $\la\in\Cbb$ and for every nontrivial solution $z(\la)$ of system \eqref{Sla}, i.e., 
 $\mL(z(\la))_k=\la\Psk\zk(\la)$ for $k\in\oinftyZ$,
  \begin{equation}\label{E:def.con}
    \sum_{k\in \mI}\zk^*(\la)\,\Psk\,\zk(\la)>0.
  \end{equation}
\end{definition}

\begin{remark}\label{R:3.2}
 Alternatively, the definiteness condition for~\eqref{Sla} can be stated in the following way: 
 System~\eqref{Sla} is definite on the discrete interval $\mI\subseteq\oinftyZ$ when, for every $\la\in\Cbb$, 
 every solution $z(\la)$ of~\eqref{Sla}, for which
  \begin{equation}\label{e3.2}
   \sum_{k\in\mI}\zk^*(\la)\,\Psk\,\zk(\la)=0,
  \end{equation}
 is trivial on $\mI$, i.e. $z_k(\la)=0,\ k\in\mI$, and as a consequence of invertibility for $\Sbb_k(\la)$ 
 (cf. \eqref{E:1A}) trivial on $\oinftyZ$. Furthermore, from the assumption that $\Psi_k\geq0, \ k\in\oinftyZ$, 
 it follows from
  \begin{equation*}\label{E:3.3}
   \sum_{k\in \mI}\zk^*(\la)\,\Psk\,\zk(\la)\leq\sum_{k\in \mhI}\zk^*(\la)\,\Psk\,\zk(\la),
   \quad\mI\subseteq\mhI\subseteq\oinftyZ,
  \end{equation*}
 so the definiteness of \eqref{Sla} on $\mI$ implies definiteness of \eqref{Sla} on every discrete interval superset
 $\mhI$, and in particular on $\oinftyZ$. Hence, definiteness of system \eqref{Sla} on some finite discrete 
 subinterval $\mI$ implies, for every $\la\in\Cbb$, that every nontrivial solution of \eqref{Sla} has a nonzero 
 semi-norm $\normP{\cdot}$. The  converse of this last statement will be shown in Lemma~\ref{L:3.3}.
\end{remark}

In the next lemma, we see that for \eqref{Sla} to be definite on a discrete interval $\mI$, it suffices to verify
\eqref{E:def.con} for nontrivial solutions of the system for only one $\la\in\Cbb$.

\begin{lemma}\label{L:3.2}
 System~\eqref{Sla} is definite on the discrete interval $\mI\subseteq\oinftyZ$ if and only if, for some $\la_0\in\Cbb$,
 each solution $z(\la_0)$ of \Sla{\la_0}, which satisfies
  \begin{equation}\label{e3.4}
   \sum_{k\in \mI}\zk^*(\la_0)\,\Psk\,\zk(\la_0)=0,
  \end{equation}
 is trivial on $\mI$, i.e., $z_k(\la_0)=0$ for $k\in\mI$.
\end{lemma}
\begin{proof}
 We begin by assuming, for $\la_0\in\Cbb$, that each solution $z(\la_0)$ of $\Sla{\la_0}$ satisfying
 \eqref{e3.4} is necessarily trivial on $\mI$, i.e., $z_k(\la_0)=0,\ k\in\mI$. Let $\la\in\Cbb$ be
 arbitrary and let $z(\la)$ be a solution of system \eqref{Sla} such that \eqref{e3.2} holds. Given
 Remark~\ref{R:3.2}, it suffices to show that $z(\la)$ is trivial on $\mI$. However, since
 $\Psk\,\zk(\la)=0$ for $k\in \mI$, we see, by \eqref{e2.9}, that $z(\la)$ is also a solution of
 system \Sla{\la_0} on $\mI$. Then, by the assumed definiteness of  \Sla{\la_0} on $\mI$,
 \eqref{e3.2} indeed implies that $z_k(\la)=0$, for $k\in\mI$. The converse is trivial.
\end{proof}

\begin{example}\label{Ex:definit}
 The particular discrete symplectic systems investigated in \cite[viz. Examples~5.1, 5.2, and 5.3]{rSH.pZ:W-TGLP} 
 are not definite, because they possess nontrivial solutions with semi-norm equal to zero. On the other hand, 
 we can provide a simple example of system \eqref{Sla} being definite on $\oinftyZ$; cf. Theorem~\ref{T:3.6B}. 
 Consider the following system:
  \begin{equation}\label{E:example}
   \mmatrix{x_k\\ u_k}=\mmatrix{1 & -1/p_{k+1}\\ \la w_k-q_k & 1+(q_k-\la w_k)/p_{k+1}}\mmatrix{x_{k+1}\\ u_{k+1}}\ 
      \text{with }\ \Psk=\mmatrix{w_k & 0\\ 0 & 0}\!,\ k\in\oinftyZ,
  \end{equation}
 where $\{p_k\}_{k=0}^\infty$, $\{q_k\}_{k=0}^\infty$, and $\{w_k\}_{k=0}^\infty$ are real-valued sequences such 
 that $w_k\geq0$, $p_{k+1}\ne0$ for $k\in\oinftyZ$ and $w_k>0$ for $k\in[a-1,b]_\sZbb$, where $a,b\in[1,\infty)_\sZbb$
 and $a\leq b$. System \eqref{E:example} corresponds to the second order Sturm--Liouville difference equation
 $\De(p_k\,\De y_{k-1})+q_k\,y_k=\la\,w_k\,y_k$. With respect to Lemma~\ref{L:3.2} it suffices to focus on the 
 definiteness of this system only for $\la=0$. Denote by $z=(x,u)^\top$ a solution of system~\eqref{E:example} with 
 $\la=0$ such that $\Psk\,z_k=0$ for $k\in\oinftyZ$, i.e., it holds
  \begin{equation}\label{E:example.id}
   \De x_k=\frac{u_{k+1}}{p_{k+1}},\qquad
   \De u_k=-q_k\,x_{k+1}+\frac{q_k}{p_{k+1}}u_{k+1},\qtextq{and}
   w_k\,x_k=0,\qquad k\in\oinftyZ.
  \end{equation}
 From the assumptions and the third condition in \eqref{E:example.id} we obtain $x_k\equiv0$ for $k\in[a-1,b]_\sZbb$ and 
 then by the first equality in \eqref{E:example.id} also $u_k\equiv0$ for $k\in[a,b]_\sZbb$. Hence $z\equiv0$
 on $[a,b]_\sZbb$, which implies $z\equiv0$ on $\oinftyZ$, because the coefficient matrix in \eqref{E:example} is 
 invertible. Therefore, system \eqref{E:example} is definite on $[a,b]_\sZbb$ and, in fact, on the whole interval 
 $\oinftyZ$ by Remark~\ref{R:3.2}.
\end{example}

For the remainder of this section, we assume that $\Phi(\la)$ represents a fundamental system of the solutions for
\eqref{Sla} such that $\Phi_{k_0}(\la)=I_{2n}$, for some $k_0\in\oinftyZ$. Note, as a consequence, that the terms 
of $\Phi(\la)$ satisfy \eqref{E:18A} for all $k\in\oinftyZ$.

The next result provides a characterization for definiteness of  \eqref{Sla} which is analogous to that for 
\eqref{E:cLHS}, as seen in \cite[Proposition~2.11]{jB.sH.hsvdS.rW11}.

\begin{lemma}\label{L:3.3}
 System~\eqref{Sla} is definite on $\oinftyZ$ if and only if there exists a finite discrete  interval
 $\mI$ over which the system is definite.
\end{lemma}
\begin{proof}
 From Remark~\ref{R:3.2} we see that definiteness of \eqref{Sla} on a finite discrete interval $\mI$ implies definiteness 
 on $\oinftyZ$. Thus, it remains to show the converse.
 
 Assume that \eqref{Sla} is definite on $\oinftyZ$.  In light of Lemma~\ref{L:3.2}, we need only to show the existence of 
 a~finite discrete interval $\mI$ over which \eqref{Sla} is definite for one value $\la_0\in\Cbb$. Thus, let 
 $\la_0\in\Cbb$ and for each finite discrete subinterval $\mI$, define the set $s(\mI)$:
  \begin{equation*}
   s(\mI)\coloneq \big\{\xi\in\Cbb^{2n}\mid\ \norm{\xi}=1,\
   \sum_{k\in \mI}\xi^*\,\Phk^*(\la_0)\,\Psk\,\Phk(\la_0)\,\xi=0\big\},
  \end{equation*}
 where $\norm{\cdot}$ denotes a norm for $\Cbb^{2n}$;  $s(\mI)$ is compact and $s({\mhI})\subseteq s(\mI)$ whenever
 $\mI\subseteq\mhI$.

 Consider a collection $\{\mI_m\}_{m\in\Nbb}$ of nested, finite, discrete intervals such that 
 $\bigcup_{m\in\Nbb} \mI_m=\oinftyZ$.  Next, assume that there exists a vector $\xi\in\Cbb^{2n}$ with $\norm{\xi}=1$ 
 such that for every $m\in\Nbb$,
  \begin{equation*}
   \sum_{k\in \mI_m}\xi^*\Phk^*(\la_0)\,\Psk\,\Phk(\la_0)\,\xi=0.
  \end{equation*}
 As a consequence,
  \begin{equation*}\label{E:3.4}
   \sum_{k\in\Nbb_0}\xi^*\Phk^*(\la_0)\,\Psk\,\Phk(\la_0)\,\xi=0.
  \end{equation*}
 Then, definiteness on $\oinftyZ$ of \eqref{Sla} implies that $\Phi_k(\la_0)\,\xi=0$ for all $k\in\oinftyZ$, and hence 
 that $\xi=0$; thus contradicting the assumption that $\|\xi\|=1$. As consequence, for the nested collection of 
 compact sets given by $\{s(\mI_m)\}_{m\in\Nbb}$,
  \begin{equation*}\label{E:3.5}
   \bigcap_{m\in\Nbb}s(\mI_m)=\emptyset.
  \end{equation*}
 Thus, $s(\mI_{m_0})=\emptyset$ for some $m_0\in\Nbb$, hence demonstrating the definiteness of \eqref{Sla} on the 
 finite discrete interval $\mI_{m_0}$.
\end{proof}

\begin{remark}
 Definiteness of \eqref{Sla} on $\oinftyZ$ plays a significant role in the development of the 
 Weyl--Titchmarsh theory for the case of a positive semi-definite weight in the definition of semi-inner 
 product. In \cite[Hypothesis~2.4]{slC.pZ10} and \cite[Hypothesis~4.11]{rSH.pZ:W-TGLP}, existence is assumed 
 for an $N_0\in\oinftyZ$ such that
  \begin{equation*}
   \sum_{k=0}^{N_0}\zk^*(\la)\,\Psk\,\zk(\la)>0
  \end{equation*}
 for every  $\la\in\Cbb$ and every nontrivial solution, $z(\la)$,  of~\eqref{Sla}. With respect to Lemma~\ref{L:3.3}, 
 this assumption is equivalent to the definiteness of \eqref{Sla} on $\oinftyZ$ for some $\la\in\Cbb$. Note also 
 that \cite[Assumption~2.2]{mB.sS10} requires satisfaction of inequality~\eqref{E:def.con} on every nonempty finite
 subinterval of $\oinftyZ$: a~condition significantly stronger than requiring definiteness of \eqref{Sla} on $\oinftyZ$ 
 as seen when \eqref{Sla} is definite on $[0,N_0]_\sZbb\subset\oinftyZ$ and  $\Psk\equiv0$ for $k\geq N_0+1$.
\end{remark}

We now give another characterization of the definiteness for~\eqref{Sla}, analogous to that given for 
systems~\eqref{E:cLHS} and \eqref{E:dLHS} in \cite[Sections~2.3 and~2.4]{mL.mmM03} and
\cite[Sections~3 and~4]{gR.yS11}, respectively.

For a discrete finite subinterval $\mI\subset\oinftyZ$, with $k_0\in\mI$, we define the $2n\times 2n$ 
positive semi-definite matrix
 \begin{equation}\label{e3.9}
  \vp(\la,\mI)\!\coloneq \sum_{k\in\mI} \Phk^*(\la)\,\Psk\, \Phk(\la).
 \end{equation}
in terms of the fundamental system $\Phi(\la)$ for \eqref{Sla} where $\Phi_{k_0}(\la)=I_{2n}$, for all $\la\in\Cbb$.  
While $\vp(\la,\mI)$ depends on $\la$ and $\mI$, we next show that the kernel and the range of $\vp(\la,\mI)$ do not 
depend on $\la$, and hence that the rank of $\vp$ is independent of $\la$.

\begin{lemma}\label{L:3.4}
 For a discrete finite subinterval $\mI\subset\oinftyZ$,  the subspaces  $\Ker \vp(\la,\mI)$ and\/
 $\Ran\vp(\la,\mI)$ are independent of $\la\in\Cbb$.
\end{lemma}
\begin{proof}
 Fix $\la\in\Cbb$ and $\xi\in\Ker\vp(\la,\mI)$, and let $z\coloneq \Phi(\la)\xi$. Then $z$ solves \Sla{\la} on $\mI$, 
 i.e., $\mL(z)_k=\la\,\Psk\,\zk$ for $k\in\mI$, while satisfying the initial condition $z_{k_0}=\xi$. As a consequence,
 $\Psk\,\zk=\Psk\,\Phk(\la)\,\xi=0$ for $k\in\mI$.  Note, for an arbitrary $\nu\in\Cbb$, that  $z$ also solves
 \Sla{\nu} on $\mI$, i.e., $\mL(z)_k=\nu\,\Psk\,\zk$, $k\in\mI$, with $z_{k_0}=\Phi_{k_0}(\nu)\xi=\xi$. Hence,
 $z=\Phi(\la)\xi=\Phi(\nu)\xi$, which implies
  \begin{equation*}
   0=\xi^*\,\vp(\la,\mI)\,\xi=\sum_{k\in\mI}\zk^*\,\Psk\,\zk=\xi^*\,\vp(\nu,\mI)\,\xi.
  \end{equation*}
 Therefore, $\Ker\vp(\la,\mI)\subseteq\Ker\vp(\nu,\mI)$. By reversing the roles of $\lambda$ and $\nu$ we obtain
 $\Ker\vp(\la,\mI)=\Ker\vp(\nu,\mI)$ for all $\la, \nu\in\Cbb$. The independence of $\Ran\vp(\la,\mI)$ on 
 $\la\in\Cbb$ follows from the fact that as defined in \eqref{e3.9}, $\vp(\la,\mI)$ is Hermitian; thus,
 $\Ran\vp(\la,\mI)=\Ran\vp^*(\la,\mI)=\ker\vp(\la,\mI)^\bot$.
\end{proof}

In general, given Lemma~\ref{L:3.4},  we shall suppress $\la$ in the following notation: 
$\Ker \vp(\la,\mI)\equiv\Ker\vp(\mI)$, $\Ran \vp(\la,\mI)\equiv\Ran \vp(\mI)$, and 
$\rank \vp(\la,\mI)\equiv\rank \vp(\mI)$ for $\la\in\Cbb$.

\begin{lemma}\label{L:5}
 There exists a discrete finite interval $\mI$, with $k_0\in\mI\subset\oinftyZ$, such that  for any discrete finite 
 interval $\mhI$ satisfying $\mI\subseteq\mhI\subset\oinftyZ$,
  \begin{equation}\label{E:3.6}
   \rank \vp(\mI)=\rank \vp(\mhI),\qquad \Ran \vp(\mI)=\Ran \vp(\mhI).
  \end{equation}
\end{lemma}
\begin{proof}
 For discrete finite intervals $\mI$ and $\mhI$ of $\oinftyZ$, where $k_0\in\mI\subseteq\mhI\subset\oinftyZ$, by 
 the definition of $\vp$ in \eqref{e3.9}, we see that $\Ker\vp(\mhI)\subseteq\Ker\vp(\mI)$. Then, given that $\vp$ 
 is Hermitian, we see that
  \begin{equation*}\label{e3.11}
   \Ran\vp(\mI) =\Ker\vp(\mI)^\bot\subseteq\Ker\vp(\mhI)^\bot=\Ran\vp(\mhI),
  \end{equation*}
 and hence that $\rank \vp(\mI)\leq \rank \vp(\mhI)$. Since $\rank\vp(\cdot)\leq 2n$, there must be a finite 
 discrete interval $\mI$ such that $ \rank \vp(\mhI)=\rank \vp(\mI)$, and $\Ran\vp(\mhI)=\Ran\vp(\mI)$ for all 
 finite discrete intervals $\mhI$ containing $\mI$.
\end{proof}

Next, we describe the connection between the definiteness of \eqref{Sla} and  the matrix $\vp(\la,\mI)$
for a finite discrete interval $\mI\subset\oinftyZ$.

\begin{theorem}\label{T:3.6}
 For a discrete finite interval $\mI\subset\oinftyZ$, and with $\vp(\la,\mI)$ defined in \eqref{e3.9}, the following
 statements are equivalent:
  \begin{itemize}
    \item[{\rm (i)}] $\rank \vp(\mI)=2n$.
    \item[{\rm (ii)}] $\Ker \vp(\mI)=\{0\}$.
    \item[{\rm (iii)}] For some $\la\in\Cbb$, every nontrivial solution, $z(\la)$, of~\eqref{Sla} satisfies
                       $\sum_{k\in\mI}\zk^*(\la)\,\Psk\,\zk(\la)>0$.
    \item[{\rm (iv)}] For some $\la\in\Cbb$,  a solution, $z(\la)$, of~\eqref{Sla} is necessarily trivial, i.e. 
                      $z_k(\la)=0$, $k\in\mI$, when $ \sum_{k\in\mI}\zk^*(\la)\,\Psk\,\zk(\la)=0$.
  \end{itemize}
\end{theorem}
\begin{proof}
 Equivalence of (i) and (ii) is clear, while equivalence of (iii) and (iv) follows from Lemma~\ref{L:3.2}. It suffices 
 then to show that statements (ii) and (iii) are equivalent.

 Assume that (ii) is true. Any nontrivial solution $z(\la)$ of~\eqref{Sla} can be expressed as $z(\la)=\Phi(\la)\,\xi$ 
 for some $\xi\in\Cbb^{2n}\setminus\{0\}$. By (ii), $\vp(\la,\mI)\,\xi\neq0$, and by the positive semi-definiteness of
 $\vp(\la,\mI)$, we see that
  \begin{equation*}
   \sum_{k\in\mI}\zk^*(\la)\,\Psk\,\zk(\la)
     =\sum_{k\in\mI} \xi^*\,\Phk^*(\la)\,\Psk\,\Phk(\la)\,\xi
     =\xi^*\,\vp(\la,\mI)\,\xi>0.
  \end{equation*}
 Conversely, assume that (iii) is true, and for $\xi\in\Cbb^{2n}\setminus\{0\}$ let $z(\la)\coloneq \Phi(\la)\,\xi$. 
 Then $z(\la)$ is a nontrivial solution of~\eqref{Sla} and, by (iii),
  \begin{equation*}
   \xi^*\,\vp(\la,\mI)\,\xi=\sum_{k\in\mI}\zk^*(\la)\,\Psk\,\zk(\la)>0.
  \end{equation*}
 Thus $\vp(\la,\mI)\,\xi\neq0$, and hence $\Ker\vp(\mI)=\{0\}$; implying that statement (ii) is satisfied.
\end{proof}

As a consequence of Lemma~\ref{L:3.3} and Theorem~\ref{T:3.6} we get the following corollary, cf.
\cite[Definition~2.14]{mL.mmM03}.

\begin{corollary}\label{C:3.6A}
 System~\eqref{Sla} is definite on $\oinftyZ$ if and only if, for some finite discrete interval $\mI\subset\oinftyZ$, 
 one of the conditions listed in Theorem~\ref{T:3.6} is satisfied.
\end{corollary}

For the special case of linear dependence on $\la$ as studied in~\cite{mB.sS10,slC.pZ10}, i.e., where the spectral 
parameter $\la$ appears only in the second equation of the system, we can show the following sufficient condition 
for the definiteness of system~\eqref{Sla} on $\oinftyZ$; q.v. Example~\ref{Ex:definit}.

\begin{theorem}\label{T:3.6B}
 Let $\la\in\Cbb$ and $\Sbb_k(\la)=\msmatrix{\A_k & \B_k\\ \C_k+\la\W_k\A_k & \D_k+\la\W_k\B_k}$ with
 $\A_k,\B_k,\C_k,\D_k,\W_k\in\Cbb^{n\times n}$ satisfy identity~\eqref{E:1.2} for all $k\in\oinftyZ$.
 If there exists an index $l\in[1,\infty)_\sZbb$ such that the matrices $\B_{l-1}$, $\W_{l-1}$, and $\W_{l}$ 
 are invertible (in fact, $\W_{l-1}$ and $\W_{l}$ are positive definite), then system~\eqref{Sla} is definite on 
 $\oinftyZ$.
\end{theorem}
\begin{proof}
 First we note that the form of the matrix $\Sbb_k(\la)$ implies $\Psk=\diag\{\W_k, 0\}$ for all $k\in\oinftyZ$
 by~\eqref{E:2}(iv). Let $z(\la)=(x(\la),\ u(\la))^\top$ be any nontrivial solution of \eqref{Sla} such that
 $\Psk\,\zk(\la)$, i.e., $\W_k\,x_k(\la)=0$, for $k\in\oinftyZ$. By Lemma~\ref{L:3.2} we have to show that $z(\la)\equiv0$ 
 on $\oinftyZ$. From the invertibility of $\W_{l-1}$ and $\W_{l}$ we obtain $x_{l-1}(\la)=x_{l}(\la)=0$. Hence
  \begin{equation*}
   0=x_{l-1}(\la)=\A_{l-1}x_{l}(\la)+\B_{l-1}u_{l}(\la)=\B_{l-1}u_{l}(\la)
  \end{equation*}
 and the invertibility of $\B_{l-1}$ implies also $u_{l}(\la)=0$, i.e., $z_{l}(\la)=0$. Since the matrix $\Sbb_k(\la)$ 
 is invertible for every $k\in\oinftyZ$, it follows $z(\la)\equiv0$ on $\oinftyZ$.
\end{proof}

%%%%%%%%%%%%%%%%%%%%%%%%%%%%%%%%%%%%%%%%%%% SECTION %%%%%%%%%%%%%%%%%%%%%%%%%%%%%%%%%%%%%%%%%%%%%%%%%%%%%%%%%%%%%%%%%%

\section{Nonhomogeneous problem}\label{S:nonhom}

The origin of the Weyl--Titchmarsh theory for discrete symplectic systems can be found in \cite{mB.sS10,slC.pZ10}, where 
the system with the spectral parameter appearing only in the second equation was studied. Recently, these results were
generalized and further extended for discrete symplectic systems with general linear dependence on the spectral parameter 
in \cite{rSH.pZ:W-TGLP}.  Given that Weyl--Titchmarsh theory has been established for discrete
symplectic systems of the form \eqref{E:1.1.1}, one can verify that these results remain valid for the
time-reversed symplectic systems \eqref{Sla} with appropriate changes in the definition of the semi-inner product and 
its weight matrix as noted in the introduction.  

In this section, we take the nonhomogeneous problem into consideration as in \cite[Section~5]{slC.pZ10} and begin by 
recalling some fundamental results from \cite{rSH.pZ:W-TGLP} which are related to the study of system \eqref{Slaf} and refer
the reader to \cite{rSH.pZ:W-TGLP} for more details.

Throughout this section, we assume that system \eqref{Sla} is definite on $\oinftyZ$, and fix the matrix $\al\in\Ga$, 
where
 \begin{equation*}
  \Ga\coloneq \big\{\al\in\Cbb^{n\times 2n}:\ \al\,\als=I,\ \al\,\J\als=0\big\}.
 \end{equation*}
By $\Ph(\la,\al)$ we denote the fundamental matrix of system \eqref{Sla} such that $\Ph_0(\la,\al)=(\als,\ -\J\als)$, 
and we emphasize its partition into $2n\times n$ blocks by the notation $\Ph(\la,\al)\coloneq (Z(\la,\al),\ \tZ(\la,\al))$. 

For $M\in\Cbb^{n\times n}$ the function
 \begin{equation}\label{E:weyl.sol}
  \Xchi_k(\la)\coloneq \Phk(\la,\al)\,\mmatrix{I & M^*}^*
 \end{equation}
represents a \emph{Weyl solution} (cf. \cite[Definition~2.11]{{rSH.pZ:W-TGLP}}), and the set
 \begin{equation*}
  D_k(\la)\coloneq \{M\in\Cbb^{n\times n}:\, \E_k(M)\leq 0\}, \qtext{where }
  \E_k(M)\coloneq i\de(\la)\,\Xchi_{k}^*(\la)\,\J\Xchi_{k}(\la)
 \end{equation*}
and $\de(\la)\coloneq\sgn\Im(\la)$, is said to be the \emph{Weyl disk}; cf. \cite[Definition~3.1]{{rSH.pZ:W-TGLP}}. 
Since system \eqref{Sla} is assumed to be definite on $\oinftyZ$, the set $D_+(\la)\coloneq \lim_{k\to\infty}D_k(\la)$ 
exists and it is closed, convex and nonempty; cf. \cite[Definition~3.10]{{rSH.pZ:W-TGLP}}. It is also well known that the 
columns of the Weyl solution $\Xchi(\la)$ defined as in \eqref{E:weyl.sol} for $k\in\oinftyZ$ through $M\in D_+(\la)$ are 
linearly independent square summable solutions of \eqref{Sla}; cf. \cite[Theorem~4.2]{{rSH.pZ:W-TGLP}}. Thus 
system~\eqref{Sla} has at least $n$ (\emph{the limit point case}) and at most $2n$ (\emph{the limit circle case}) linearly 
independent square summable solutions. 

By $M_+(\la)$ we denote the \emph{half-line Weyl--Titchmarsh $M(\la)$-functions}, defined in accordance with 
\cite[Remark~3.17]{rSH.pZ:W-TGLP}, and note that
 \begin{equation}\label{E:mpla}
  M_+^*(\la)=M_+(\bla),\qquad \la\in\Cbb\setminus\Rbb.
 \end{equation}
If systems \eqref{Sla} and \Sla{\nu}, where $\la,\nu\in\Cbb\setminus\Rbb$, are both in the limit point or in the limit
circle case, then
 \begin{equation}\label{E:weyl.limit}
   \lim_{k\to\infty}\Xchi_k^{+*}(\la)\,\J\Xchi_k^{+}(\nu)=0,
 \end{equation}
where $\Xchi^{+}(\la)$ and $\Xchi^{+}(\nu)$ represent Weyl solutions for systems \eqref{Sla} and \Sla{\nu} defined 
by \eqref{E:weyl.sol} for $k\in\oinftyZ$, corresponding to the matrices $M_+(\la)$ and $M_+(\nu)$, respectively; cf.
\cite[Theorem~4.12]{rSH.pZ:W-TGLP}. 

For $\la\in\Cbb\setminus\Rbb$ and $k,l\in\oinftyZ$ we introduce the \emph{Green function}
 \begin{equation}\label{E:green.def.kl}
   G_{k,l}(\la)\coloneq \begin{cases}
                   \tZ_k(\la)\,\Xchi_l^{+*}(\bla), &k\in[0,l]_\sZbb,\\
                   \Xchi^+_k(\la)\,\tZ_l^*(\bla),  &k\in[l+1,\infty)_\sZbb.
                 \end{cases}
 \end{equation}
In the literature, we also find the terminology \emph{resolvent kernel} for an analogous function in the continuous time 
case; cf. \cite[page~15]{viK.fsR76}. The function $G_{k,l}(\la)$ can be equivalently written as
 \begin{equation}\label{E:green.def.lk}
   G_{k,l}(\la)=\begin{cases}
                   \Xchi^+_k(\la)\,\tZ_l^*(\bla),  &l\in[0,k-1]_\sZbb,\\
                   \tZ_k(\la)\,\Xchi_l^{+*}(\bla), &l\in[k,\infty)_\sZbb.
                   \end{cases}
 \end{equation}

\begin{lemma}
 Let $\al\in\Ga$, $\la\in\Cbb\setminus\Rbb$, and system \eqref{Sla} be definite on $\oinftyZ$. Then
  \begin{equation}\label{E:5.50}
    \Xchi^+_k(\la)\,\tZ_k^*(\bla)-\tZ_k(\la)\,\Xchi_k^{+*}(\bla)=\J \qtext{for all} k\in\oinftyZ.
  \end{equation}
\end{lemma}
\begin{proof}
 Identity \eqref{E:5.50} then follows by a direct calculation from the definition of $\Xchi^+(\cdot)$ and identities
 \eqref{E:18A}(iii) and \eqref{E:mpla}.
\end{proof}

In the next lemma, some fundamental properties of the Green function, $G(\la)$, are established. We note that the 
given identities are presented in a more symmetric form, with respect to the variables $k,l$, than the corresponding
identities for the Green function in the case of the discrete symplectic system with the special linear dependence 
on the spectral parameter given in \cite[Lemma~5.1]{slC.pZ10}.

\begin{lemma}\label{L:green.prop}
 Let $\al\in\Ga$, $\la\in\Cbb\setminus\Rbb$, and system \eqref{Sla} be definite on $\oinftyZ$. Then the function
 $G_{\cdot,\cdot}(\la)$ possesses the following properties:
  \begin{itemize}
    \item[\rm(i)] $G_{k,l}^*(\la)=G_{l,k}(\bla)$ for all $k,l\in\oinftyZ$ such that $k\ne l$;
    \item[\rm(ii)] $G_{k,k}^*(\la)=G_{k,k}(\bla)+\J$ for all $k\in\oinftyZ$;
    \item[\rm(iii)] for every $k,l\in\oinftyZ$ such that $k\in\T(l)$, the function $G_{\cdot,l}(\la)$ solves the 
                    homogeneous system \eqref{Sla} on the set $\T(l)$, where
                     \begin{equation*}
                      \T(l)\coloneq \{\tau\in\oinftyZ:\ \tau\ne l\};
                     \end{equation*}
    \item[\rm(iv)] $G_{k,k}(\la)=\Sbb_k(\la)\,G_{k+1,k}(\la)-\J$ for every $k\in\oinftyZ$;
    \item[\rm(v)] the columns of $G_{\cdot,l}(\la)$ belong to $\ellP$ for every $l\in\oinftyZ$ and the columns of
                  $G_{k,\cdot}(\la)$ belong to $\ellP$ for every $k\in\oinftyZ$.
  \end{itemize}
\end{lemma}
\begin{proof}
 The first property follows directly from the definition of $G_{k,l}(\la)$ in \eqref{E:green.def.kl}. The second property 
 can be obtained from \eqref{E:green.def.kl} by means of identity \eqref{E:5.50}. To prove the third property, it is
 necessary to distinguish between $k$ and $l$ in the relation: the statement follows from the fact that the functions
 $\Xchi^+(\la)$ and $\tZ(\la)$ solve system \eqref{Sla} with respect to $k$. Property (iv) can be proven using the
 definition of $G(\la)$ in \eqref{E:green.def.kl} and identities \eqref{E:mpla} and \eqref{E:18A}. Finally, the columns of
 $G_{\cdot,l}(\la)$ belong to $\ellP$ for every $l\in\oinftyZ$ by the definition of the Green function, because
  \begin{align*}
   \normP{G_{\cdot,l}(\la)\,e_j}^2=
      &\,e_j^*\,\Xchi^+_l(\bla)\,\bigg(\sum_{k=0}^l\tZ_k^*(\la)\,\Psk\,\tZ_k(\la)\bigg)\,\Xchi_l^{+*}(\bla)\,e_j\\
      &+e_j^*\,\tZ_l(\bla)\,\bigg(\sum_{k=l+1}^\infty\Xchi_k^{*+}(\la)\,\Psk\,\Xchi^+_k(\la)\bigg)\tZ_k^*(\bla)\,e_j<\infty,
  \end{align*}
 while the columns of $G_{k,\cdot}(\la)$ are in $\ellP$ for every $k\in\oinftyZ$ by a similar calculation 
 and~\eqref{E:green.def.lk}. 
\end{proof}

We associate with the nonhomogeneous system \eqref{Slaf} the function
 \begin{equation}\label{E:5.54}
   \hz_k(\la)=\sum_{l=0}^\infty G_{k,l}(\la)\,\Ps_l\,f_l,\qquad k\in\oinftyZ,
 \end{equation}
which is well defined for all $f\in\ellP$, because the columns of $G_{k,\cdot}(\la)$ are square summable by the previous 
lemma. In addition, by \eqref{E:green.def.lk}, we can write
 \begin{equation}\label{E:5.55}
  \hz_k(\la)=\bigg\{\sum_{l=0}^{k-1}+\sum_{l=k}^{\infty}\bigg\}\,G_{k,l}(\la)\,\Ps_l\,f_l
   =\Xchi^{+}_k(\la)\sum_{l=0}^{k-1}\tZ^*_l(\bla)\,\Ps_l\,f_l+\tZ_k(\la)\sum_{l=k}^{\infty}\Xchi_l^{+*}(\bla)\,\Ps_l\,f_l.
 \end{equation}
Similarly as in \cite[Theorem~5.2]{slC.pZ10} we show that the above defined function $\hz(\la)$ represents a square summable
solution of system \eqref{Slaf}.

\begin{theorem}\label{T:5.3}
 Let $\al\in\Ga$, $\la\in\Cbb\setminus\Rbb$, $f\in\ellP$, and system \eqref{Sla} be definite on $\oinftyZ$. The function
 $\hz(\la)$ defined in \eqref{E:5.54} solves system \eqref{Slaf}, satisfies the initial condition $\al\,z_0(\la)=0$, is
 square summable, i.e., $\hz(\la)\in\ellP$, and it holds
  \begin{equation}\label{E:5.56A}
   \normP{\hz(\la)}\leq\frac{1}{\abs{\Im(\la)}}\normP{f}.
  \end{equation}
 In addition, if system \eqref{Sla}  is in the limit point or limit circle case for all $\la\in\Cbb\setminus\Rbb$, then
  \begin{equation}\label{E:5.56B}
   \lim_{k\to\infty} \Xchi_k^{+*}(\nu)\,\J\,\hz_k(\la)=0 \qtext{for every} \nu\in\Cbb\setminus\Rbb.
  \end{equation}
\end{theorem}
\begin{proof}
 The form of $\hz_k(\la)$ given in \eqref{E:5.55} together with a similar expression of $\hz_{k+1}(\la)$, the facts that
 $\Xchi^+(\la)$ and $\tZ(\la)$ solve \eqref{Sla}, and identity~\eqref{E:5.50} yield
  \begin{align*}
   \hz_k(\la)-\Sbb_k(\la)\,\hz_{k+1}(\la)
      &=\Xchi^+_k(\la)\sum_{l=0}^{k-1}\tZ^*_l(\bla)\,\Ps_l\,f_l+\tZ_k(\la)\sum_{l=k}^{\infty}\Xchi_l^{+*}(\bla)\,\Ps_l\,f_l
      \\[-2mm]
       &\hspace*{5mm}-\Sbb_k(\la)\,\Xchi^+_{k+1}(\la)\sum_{l=0}^{k}\tZ^*_l(\bla)\,\Ps_l\,f_l
       -\Sbb_k(\la)\,\tZ_{k+1}(\la)\sum_{l=k+1}^{\infty}\Xchi_l^{+*}(\bla)\,\Ps_l\,f_l\\
      &=-\big[\Xchi^+_k(\la)\,\tZ_k^*(\bla)-\tZ_k(\la)\,\Xchi_k^{+*}(\bla)\big]\Psk\,f_k=-\J\,\Psk\,f_k,
  \end{align*}
 i.e., the function $\hz(\la)$ solves system \eqref{Slaf}. The fulfillment of the boundary condition follows by the simple
 calculation
  \begin{equation*}
   \al\,\hz_0(\la)=\al\,\tZ_0(\la)\sum_{l=0}^\infty\Xchi_l^{+*}(\bla)\,\Ps_l\,f_l
                  =-\al\,\J\als\sum_{l=0}^\infty\Xchi_l^{+*}(\bla)\,\Ps_l\,f_l
                  =0,
  \end{equation*}
 because $\tZ_0(\la)=-\J\als$ and $\al\in\Ga$. 

 Next, we prove the estimate in \eqref{E:5.56A} which together with the assumption $f\in\ellP$ will imply that
 $\hz(\la)\in\ellP$. For every $r\in\oinftyZ$ we define the function 
  \begin{equation*}
   f_k^{[r]}\coloneq \begin{cases}
               f_k, &k\in[0,r]_\sZbb,\\
               0,   &k\in[r+1,\infty)_\sZbb,
              \end{cases}
  \end{equation*}
 and the function
  \begin{equation*}
   \hz_k^{[r]}(\la)\coloneq \sum_{l=0}^{\infty}G_{k,l}(\la)\,\Ps_l\,f_l=\sum_{l=0}^{r}G_{k,l}(\la)\,\Ps_l\,f_l.
  \end{equation*}
 The function $\hz^{[r]}_k(\la)$ solves system \eqref{Slaf} with $f$ replaced by 
 $f^{[r]}$. Applying the extended Lagrange identity from Theorem~\ref{T:1}, we obtain
  \begin{equation}\label{E:5.58A}
   \begin{aligned}
    \lim_{k\to\infty}\hz_{k+1}^{[r]*}(\la)\,\J\,\hz_{k+1}^{[r]}(\la)
       &=\hz_{0}^{[r]*}(\la)\,\J\,\hz_{0}^{[r]*}(\la)+(\bla-\la)\sum_{k=0}^{\infty}\hz_k^{[r]*}(\la)\,\Psk\,\hz_k^{[r]}(\la)\\    
       &\hspace*{15mm}+\sum_{k=0}^{\infty}f_k^{[r]*}\,\Psk\,\hz_k^{[r]}(\la)
                      -\sum_{k=0}^{\infty}\hz_k^{[r]*}(\la)\,\Psk\,f_k^{[r]}.
   \end{aligned}
  \end{equation}
 Since $\tZ_0(\la)=-\J\als$ and $\al\in\Ga$, we see that
  \begin{equation*}
   \hz_{0}^{[r]*}(\la)\,\J\,\hz_{0}^{[r]}(\la)
    =\bigg(\sum_{l=0}^{r}\Xchi_l^{+*}(\bla)\,\Ps_l\,f_l\bigg)^{\!\!*}\tZ_0^*(\la)\,\J\,\tZ_0(\la)
                                                     \bigg(\sum_{l=0}^{r}\Xchi_l^{+*}(\bla)\,\Ps_l\,f_l\bigg)=0,
  \end{equation*}
 where, for every $k\in[r+1,\infty)_\sZbb$, we can write
  \begin{equation}\label{E:5.58B}
   \hz^{[r]}_k(\la)=\Xchi^+_k(\la)\,g_r(\la),\qtext{where} g_r(\la)\coloneq \sum_{l=0}^{r}\tZ_{l}^*(\bla)\,\Ps_l\,f_l.
  \end{equation}
 This, together with the fact $M_+(\la)\in D_+(\la)$, yields
  \begin{align*}
    \frac{1}{\bla-\la}\,\lim_{k\to\infty}\hz_{k+1}^{[r]*}(\la)\,\J\,\hz_{k+1}^{[r]}(\la)
      &=\frac{i\de(\la)}{2\abs{\Im(\la)}}\,g_r^*(\la)
                        \bigg(\lim_{k\to\infty}\Xchi_{k+1}^{+*}(\la)\,\J\,\Xchi^+_{k+1}(\la)\bigg)\,g_r(\la)\\
      &=\frac{1}{2\abs{\Im(\la)}}\,g_r^*(\la)\Big(\lim_{k\to\infty}\E_{k+1}(M_+(\la))\Big)\,g_r(\la)\leq0.
  \end{align*}
 Using the Cauchy--Schwarz inequality, $\Ps\geq0$, and  identity \eqref{E:5.58A}, we see that
  \begin{align*}
   \normP{\hz^{[r]}(\la)}^2=\sum_{k=0}^{\infty}\hz_k^{[r]*}(\la)\,\Psk\,\hz_k^{[r]}(\la)
      &\leq \frac{1}{2i\Im(\la)}\bigg(\sum_{k=0}^{r}f_k^{[r]*}\,\Psk\,\hz_k^{[r]}(\la)
                                      -\sum_{k=0}^{r}\hz_k^{[r]*}(\la)\,\Psk\,f_k^{[r]}\bigg)\\
      &\leq\frac{1}{\abs{\Im(\la)}}\bigg|\sum_{k=0}^r \hz_k^{[r]*}(\la)\,\Psk\,f_k^{[r]}\bigg|\\
      &\leq\frac{1}{\abs{\Im(\la)}}\bigg(\sum_{k=0}^r \hz_k^{[r]*}(\la)\,\Psk\,\hz_k^{[r]}(\la)\bigg)^{\!\!1/2}\,
                                   \bigg(\sum_{k=0}^r f_k^{[r]*}\Psk\,f_k^{[r]}\bigg)^{\!\!1/2}\\
      &\leq\frac{1}{\abs{\Im(\la)}}\normP{\hz^{[r]}(\la)}\,\normP{f^{[r]}},
  \end{align*}
 thereby yielding the inequality
  \begin{equation}\label{E:5.59}
   \normP{\hz^{[r]}(\la)}\leq\frac{1}{\abs{\Im(\la)}}\normP{f^{[r]}}\leq\frac{1}{\abs{\Im(\la)}}\normP{f}.
  \end{equation} 
 For any $k,r\in\oinftyZ$, we now see that
  \begin{equation*}
   \hz_k(\la)-\hz_k^{[r]}(\la)=\sum_{l=r+1}^\infty G_{k,l}(\la)\,\Ps_l\,f_l.
  \end{equation*}

 Let $m\in [0,r]_\sZbb$ be fixed. By the definition of $G(\la)$ in \eqref{E:green.def.kl}, we obtain for every 
 $k\in[0,m]_\sZbb$ that
  \begin{equation}\label{E:5.60}
   \hz_k(\la)-\hz_k^{[r]}(\la)=\tZ_k(\la)\sum_{l=r+1}^\infty\Xchi_l^{+*}(\bla)\,\Ps_l\,f_l.
  \end{equation}
 Since the columns of $\Xchi^+(\bla)$ and the function $f$ belong to $\ellP$, it follows that the right-hand side of
 \eqref{E:5.60} tends to zero as $r\to\infty$ for every $k\in[0,m]_\sZbb$. Hence $\hz^{[r]}$ converges uniformly to the
 function $\hz(\la)$ on the interval $[0,m]_\sZbb$. %Since $\hz(\la)=\hz^{[r]}(\la)$ on $[0,m]_\sZbb$ and $\Ps\geq0$, 
 We see by \eqref{E:5.59} that
  \begin{equation*}
   \sum_{k=0}^m \hz_k^{[r]*}(\la)\,\Psk\,\hz_k^{[r]}(\la)\leq\normP{\hz^{[r]}(\la)}^2\leq\frac{1}{\abs{\Im(\la)}^2}\,\normP{f}^2.
  \end{equation*}
 Then, as a consequence of the uniform convergence for $r\to\infty$ on $[0,m]_\sZbb$, we see that
  \begin{equation}\label{E:5.61}
   \sum_{k=0}^m \hz_k^*(\la)\,\Psk\,\hz_k(\la)\leq\frac{1}{\abs{\Im(\la)}^2}\,\normP{f}^2.
  \end{equation}
 Upon taking the limit for $m\to\infty$ in \eqref{E:5.61} the desired estimate in \eqref{E:5.56A} follows. 
 
 Finally, to establish the existence of the limit in \eqref{E:5.56B}, assume that system \eqref{Sla} is in the limit point
 case for all $\la\in\Cbb\setminus\Rbb$. From the extended Lagrange identity in Theorem~\ref{T:1}, for any $k,r\in\oinftyZ$, 
 we obtain
  \begin{equation}\label{E:5.61A}
   \Big[\Xchi_{j}^{+*}(\nu)\,\J\,\hz^{[r]}_{j}(\la)\Big]_{0}^{k+1}
    =(\bar{\nu}-\la)\sum_{j=0}^{k}\Xchi_{j}^{+*}(\nu)\,\Ps_j\,\hz^{[r]}_{j}(\la)
              -\sum_{j=0}^{k}\Xchi_{j}^{+*}(\nu)\,\Ps_j\,f^{[r]}_{j}(\la).
  \end{equation}
 For $k\in[r+1,\infty)_\sZbb$ identity \eqref{E:5.58B} holds, and,  by \eqref{E:weyl.limit}, we obtain
  \begin{equation*}
   \lim_{k\to\infty}\Xchi_{k+1}^{+*}(\nu)\,\J\,\hz^{[r]}_{k+1}(\la)
    =\lim_{k\to\infty}\Xchi_{k+1}^{+*}(\nu)\,\J\,\Xchi^+_{k+1}(\la)\,g_r(\la)=0.
  \end{equation*}
 In the limit, as $k\to\infty$,  \eqref{E:5.61A}  yields
  \begin{equation}\label{E:5.63}
   \Xchi_{0}^{+*}(\nu)\,\J\,\hz^{[r]}_{0}(\la)
     =(\la-\bar{\nu})\sum_{j=0}^{\infty}\Xchi_{j}^{+*}(\nu)\,\Ps_j\,\hz^{[r]}_{j}(\la)
       +\sum_{j=0}^{\infty}\Xchi_{j}^{+*}(\nu)\,\Ps_j\,f^{[r]}_{j}(\la).
  \end{equation}
 It was proven earlier that $\hz^{[r]}(\la)$ converges uniformly  on finite subintervals of $\oinftyZ$  as $r\to\infty$. 
 Then, from \eqref{E:5.63}, we see that 
  \begin{equation}\label{E:5.63A}
   \Xchi_{0}^{+*}(\nu)\,\J\,\hz_{0}(\la)
     =(\la-\bar{\nu})\sum_{j=0}^{\infty}\Xchi_{j}^{+*}(\nu)\,\Ps_j\,\hz_{j}(\la)
       +\sum_{j=0}^{\infty}\Xchi_{j}^{+*}(\nu)\,\Ps_j\,f_{j}(\la).
  \end{equation}
 On the other hand, by the extended Lagrange identity, for every $k\in\oinftyZ$,  we obtain
  \begin{equation}\label{E:5.64}
     \Big[\Xchi_{j}^{+*}(\nu)\,\J\,\hz_{j}(\la)\Big]_{0}^{k+1}
        =(\bar{\nu}-\la)\sum_{j=0}^{k}\Xchi_{j}^{+*}(\nu)\,\Ps_j\,\hz_{j}(\la)
          -\sum_{j=0}^{k}\Xchi_{j}^{+*}(\nu)\,\Ps_j\,f_{j}(\la).
  \end{equation}
 Letting $k\to\infty$ in \eqref{E:5.64}, and using inequality \eqref{E:5.63A}, the limit in \eqref{E:5.56B} is established. 
 
 An argument, similar to that given above, can be used in the limit circle case to show the existence of the limit in
 \eqref{E:5.56B}, because all solutions of \eqref{Sla} are square summable. However, an alternative and more 
 direct method of the proof is available. 

 By \eqref{E:5.55} we have for every $k\in\oinftyZ$ that
  \begin{equation}\label{E:5.65}
   \Xchi_k^{+*}(\nu)\,\J\,\hz_k(\la)=
     \Xchi_k^{+*}(\nu)\,\J\,\Xchi_k^+(\la)\sum_{l=0}^{k-1}\tZ^*_l(\bla)\,\Ps_l\,f_l
     +\Xchi_k^{+*}(\nu)\,\J\,\tZ_k(\la)\sum_{l=k}^{\infty}\Xchi_l^{+*}(\bla)\,\Ps_l\,f_l.
  \end{equation}
 The limit of the first term on the right-hand side in \eqref{E:5.65} is equal to zero because
 $\Xchi_k^{+*}(\nu)\,\J\,\Xchi_k^+(\la)$ tends to zero by \eqref{E:weyl.limit} and it is multiplied by a sum which
 converges as $k\to\infty$. The second term on the right-hand side of \eqref{E:5.65} tends also to zero, because the
 columns of $\tZ(\la)$ are square summable. This implies that the function $\Xchi^{+*}(\nu)\,\J\,\tZ(\la)$ is bounded 
 and it is multiplied by a sum converging to zero as $k\to\infty$; thereby establishing the limit in \eqref{E:5.56B}.
\end{proof}

In the last result of this section, we extend \cite[Corollary~5.3]{slC.pZ10} to the case of general linear dependence 
on the spectral parameter.

\begin{corollary}\label{C:5.4}
 Let $\al\in\Ga$, $\la\in\Cbb\setminus\Rbb$, $f\in\ellP$, and $v\in\Cbb^n$. Let system \eqref{Sla} be definite on
 $\oinftyZ$, and define the function
  \begin{equation}\label{E:5.70}
   \hy_k(\la)\coloneq \Xchi_k^+(\la)\,v+\hz_k(\la),\quad k\in\oinftyZ,
  \end{equation}
 where $\hz(\la)$ is given in \eqref{E:5.54}. Then, $\hy(\la)$ represents a square summable solution of system \eqref{Slaf}
 satisfying the initial condition $\al\,\hy_0(\la)=v$ and 
  \begin{equation}\label{E:5.71}
   \normP{\hy(\la)}\leq \frac{1}{\abs{\Im(\la)}}\,\normP{f}+\normP{\Xchi^+(\la)\,v}.
  \end{equation}
 If system \eqref{Sla} is in the limit point or in the limit circle case for all $\la\in\Cbb\setminus\Rbb$, we have
  \begin{equation}\label{E:5.72}
   \lim_{k\to\infty} \Xchi_k^{+*}(\nu)\,\J\,\hy_k(\la)=0 \qtext{for every} \nu\in\Cbb\setminus\Rbb.
  \end{equation}
 Moreover, in the limit point case, the function $\hy(\la)$ is the unique square summable solution of system 
 \eqref{Slaf} satisfying $\al\,\hy_0(\la)=v$, while in the limit circle case $\hy(\la)$ is the unique solution of
 \eqref{Slaf} being in $\ellP$ such that $\al\,\hy_0(\la)=v$ and
  \begin{equation}\label{E:5.73}
   \lim_{k\to\infty}\Xchi^{+*}_k(\bla)\,\J\,\hy_k(\la)=0.
  \end{equation}
\end{corollary}
\begin{proof}
 Since $\Xchi^+(\la)\,v$ solves system \eqref{Sla} and $\Ph_0(\la,\al)=(\als,\ -\J\als)$, it follows from 
 Theorem~\ref{T:5.3} that $\hy(\la)$ solves the nonhomogeneous system \eqref{Slaf} and satisfies
 $\al\,\hy_0(\la)=\al\,\Xchi_0^+(\la)\,v=v$. The estimate in \eqref{E:5.71} follows directly from \eqref{E:5.70} 
 and \eqref{E:5.56A} by the triangle inequality. The limit in \eqref{E:5.72} follows in the limit circle or in the 
 limit point case from \eqref{E:weyl.limit}, \eqref{E:5.56B}, and from the calculation
  \begin{equation*}
   \lim_{k\to\infty}\Xchi_k^{+*}(\nu)\,\J\,\hy_k(\la)
     =\lim_{k\to\infty}\big\{\Xchi_k^{+*}(\nu)\,\J\,\Xchi_k^+(\la)\,v+\Xchi_k^{+*}(\nu)\,\J\,\hz_k(\la)\big\}=0.
  \end{equation*}
 Finally, we prove the uniqueness of the solution in the limit point and in the limit circle case. Assume that 
 $y^{[1]}(\la)$ and $y^{[2]}(\la)$ are two square summable solutions of \eqref{Slaf} satisfying
 $\al\,y^{[1]}_0(\la)=v=\al\,y^{[2]}_0(\la)$. Then the function $y_k(\la)\coloneq y^{[1]}_k(\la)-y^{[2]}_k(\la)$,
 $k\in\oinftyZ$, represents a square summable solution of system \eqref{Sla}, which satisfies $\al\,y_0(\la)=0$. 
 Since $y(\la)=\Ph(\la,\al)\,c$ for some $c\in\Cbb^{2n}$, the initial condition $\al\,y_0(\la)=0$
 implies that $y_k(\la)=\tZ_k(\la)\,d$ for some $d\in\Cbb^n$. If system \eqref{Sla} is in the limit point case, we 
 have $y(\la)\not\in\ellP$ for $d\ne0$, because the columns of $\tZ(\la)$ do not belong to $\ellP$ in this case; 
 cf. \cite[Theorem~4.4]{rSH.pZ:W-TGLP}. Therefore $d=0$ and the uniqueness follows. On the other hand, in the limit 
 circle case we obtain from the previous part of the proof, from the limit in \eqref{E:5.73}, and from the first 
 identity in \eqref{E:18A} that
  \begin{equation*}
   0=\lim_{k\to\infty}\Xchi^{+*}_k(\bla)\,\J\,\Ph_k(\la,\al)\mmatrix{0\\ d}
    =\lim_{k\to\infty}\mmatrix{I & M_+^*(\bla)}\Ph_k^*(\bla,\al)\,\J\,\Ph_k(\la,\al)\mmatrix{0\\ d}
    =d,
  \end{equation*}
 which implies the uniqueness of the solution $\hy(\la)$ also in this case.
\end{proof}

%%%%%%%%%%%%%%%%%%%%%%%%%%%%%%%%%%%%%%%%%%% SECTION %%%%%%%%%%%%%%%%%%%%%%%%%%%%%%%%%%%%%%%%%%%%%%%%%%%%%%%%%%%%%%%%%%

\section{Maximal and minimal linear relations}\label{S:4}
\subsection{Linear relations and definiteness}
We now return to the topic of linear relations introduced in Section~\ref{S:2.2} and focus on a pair of linear relations 
defined in terms of the linear map, $\mL$, introduced in \eqref{e2.9} in association with the discrete symplectic system
\eqref{Sla}. We point out that similar results in association with linear Hamiltonian differential and difference systems 
can be found in \cite{mL.mmM03,jB.sH.hsvdS.rW11,gR.yS11}. We first introduce the maximal linear relation, $\Tmax$, as a
subspace of $\tellP^2\coloneq\tellP\times\tellP$ defined by
 \begin{equation}\label{E:4.1}
  \Tmax\coloneq \big\{\{\tz,\tf\}\in\tellP^2\mid \exists\ u\in \tz \ \text{such that} \ \mL(u)=\Psi f \big\}.
 \end{equation}
Note that when $\mL(u)=\Psi f$, then $\mL(u)=\Psi g$ for all $g\in\tf$. Similarly, we define a  pre-minimal linear
relation, $\Tnod$, by
 \begin{equation}\label{E:4.2}
  \Tnod\coloneq \big\{\{\tz,\tf\}\in\tellP^2\mid\exists\ u\in\tz\cap\ellPn\ \text{such that} \ 
                    \mL(u)=\Psi f \big\}\subseteq \Tmax,
 \end{equation}
and, by \eqref{e2.18a}, we define
 \begin{equation}\label{e5.3}
  \Tnod-\lambda I\coloneq \big\{\{\tz,\tf\}\in\tellP^2\mid\exists\ u\in\tz\cap\ellPn\ \text{such that} \
  \mL(u)-\la\Psi u=\Psi f \big\}\subseteq\Tmax -\la I.
 \end{equation}

The consideration of linear relations in our current context is natural given that the weight represented by $\Psi$, 
present in \eqref{Slaf} and the sequence spaces associated with $\Tmax$ and $\Tmin$, has terms none of which are 
positive definite, but all of which are positive semi-definite; cf. \eqref{e2.2}. A simple example for \eqref{Slaf}, 
analogous to that found in  \cite[Section~2]{mL.mmM03} in association with \eqref{E:cLHS}, illustrates.  

\begin{example}\label{E:non-den.op}
 For the system \Slaf{0}{f} with $n=1$, let 
  \begin{equation}\label{e5.4}
   \S_k=\begin{pmatrix} 1&0\\0&1\end{pmatrix},\quad \Psk=\begin{pmatrix} 1&0\\0&0\end{pmatrix},
   \quad f_k=\begin{pmatrix}f_k^{[1]}\\f_k^{[2]}\end{pmatrix},\quad
   z_k=\begin{pmatrix}x_k\\u_k\end{pmatrix},\quad k\in\Nbb_0.
  \end{equation}
 Then $\mL(z)$ and \Slaf{0}{f}, respectively, can be written as
  \begin{equation}\label{E:example.oper}
   \mL(z)=\begin{pmatrix} 0&-1\\1&\phantom{-}0 \end{pmatrix}  \De z, \qquad
   \De z= \begin{pmatrix} \De x\\ \De u\end{pmatrix}=\begin{pmatrix}0\\-f^{[1]}\end{pmatrix}.
  \end{equation}
 Then, for any $f\in\ellP$, 
  \begin{equation*}\label{e5.5}
   z_0=0,\quad z_k=\begin{pmatrix} 0\\ -\sum_{j=0}^{k-1}f_j^{[1]}\end{pmatrix},\quad k\in\Nbb, 
  \end{equation*}
 is a solution of \Slaf{0}{f} in \eqref{E:example.oper}. In particular, $z\in\tilde 0$ and as a consequence, 
 $\{\tilde0,\tf\}\in\Tmax$ for any $\tf\ne\tilde 0$. Hence, $\mul\Tmax$, as defined in Section~\ref{S:2.2}, is 
 nontrivial, and $\Tmax$, while a linear relation, is not a~linear operator.  Note also that a solution of \Slaf{0}{f} 
 in \eqref{E:example.oper} which is an element of $\ellPn$ of necessity is such that $x_k=0$ for all $k\in\Nbb_0$ with 
 the consequence that $\dom T_0 =\{\tilde 0\}$.
\end{example}

The system \eqref{Sla} given by \eqref{e5.4} in Example~\ref{E:non-den.op} is not definite. The next result, 
characterizes definiteness of \eqref{Sla} in terms of the domain of $\Tmax$.

\begin{theorem}
 System~\eqref{Sla} is definite on $\oinftyZ$ if and only if for any $\{\tz,\tf\}\in\Tmax$ there exists a~unique 
 $u\in\tz$ such that $\mL(u)=\Psi f$.
\end{theorem}
\begin{proof}
 Assume that system \eqref{Sla} is definite. Let $\{\tz,\tf\}\in\Tmax$, let $z^{[1]}, z^{[2]}\in\tz$, and let 
 $y\coloneq z^{[1]}- z^{[2]}$. Then $y\in\tilde 0\in\tellP$, and hence $\Psk y_k=0$, for all $k\in\Nbb_0$. 
 Moreover, $\mL(y)=0$ implies that $y$ is a~solution of \eqref{Sla} for $\la=0$. Then, by Remark~\ref{R:3.2}, 
 the definiteness of \eqref{Sla} implies that $y=0$.

 To show the converse, assume that there is only one $u\in\tz$ for which $\mL(u)=\Psi f$, when $\{\tz,\tf\}\in\Tmax$. 
 Let $\mI\subset\Nbb_0$ be a discrete finite interval such that $\rank\vp(\mI)$ is maximal; cf. Lemma~\ref{L:5}. 
 If $\rank\vp(\mI)<2n$, then there is $\eta\in\Cbb^{2n}\setminus\{0\}$ such that
 $\eta^*\vp(0,\mI)\eta=\sum_{k\in\mI}\eta^*\Phk^*\Psk\Phk\eta=0$, where $\Phk=\Phi_k(0)$. If
 $\sum_{k\in\Nbb_0}\eta^*\Phk^*\Psk\Phk\eta=0$, then $u=\Ph\eta\in\ellP$. Moreover, $\mL(u)=0$, and $u\in\tilde 0$. 
 Given that $0$ is the unique sequence in $\tilde 0$ satisfying $\mL(z)=0$, then $u=\Ph\eta=0$ and as a consequence,
 $\eta=0$, which contradicts the assumption that $\eta\ne 0$. Hence, there is a discrete finite interval superset, $\mhI$, 
 of $\mI$ such that $\vp(0,\mhI)\eta\ne 0$. As a result, $\ker\vp(\mhI)\subset\ker\vp(\mI)$, and hence,
 $\ran\vp(\mI)\subset\ran\vp(\mhI)$; thereby contradicting the maximality of $\rank\vp(\mI)$. Thus, $\rank\vp(\mI)=2n$ and
 \eqref{Sla} is definite by Theorem~\ref{T:3.6}.
\end{proof}

\begin{remark}
Example~\ref{E:non-den.op} shows that the linear relation given by $\mL$ does not determine an operator; neither its 
domain is dense in $\tellP$. Moreover, consider the definite system in \eqref{E:example} from Example~\ref{Ex:definit} 
with $p_k\equiv1$, $q_k\equiv0$, and $w_k\equiv1$ for all $k\in\oinftyZ$. If we put $f_0=\msmatrix{1\\0}$ and
$f_k=\msmatrix{0\\0}$ for $k\in\Nbb$ with the notation $f_k=\big(f^{[1]}_k,\ f^{[2]}_k\big)^{\!\top}$ for 
$k\in\oinftyZ$, then the corresponding nonhomogeneous system \Slaf{0}{f}, i.e.,
 \begin{equation*}
  \mmatrix{x_k\\ u_k}=\mmatrix{1 & -1\\ 0 & \phantom{}1}\mmatrix{x_{k+1}\\ u_{k+1}}+\mmatrix{0\\ f^{[1]}_k},
 \end{equation*}
possesses the solution $z_k=\msmatrix{x_k\\ u_k}$ such that $z_0=\msmatrix{0\\ 1}$ and $z_k=\msmatrix{0\\0}$ for 
$k\in\Nbb$, i.e.,  there exists $\tf\neq\tilde0$ such that $\{\tilde0,\tf\}\in\Tmax$. This shows that 
definiteness of system \eqref{Sla} does not suffice to prove that $\mul\Tmax=\{\tilde0\}$, cf. \cite{yS.hS11}. 

Thus, to guarantee that the linear relation defines an operator, we need to assume explicitly that 
$\mul\Tmax=\{\tilde0\}$, i.e., if there exists $z\in\tilde{0}$ such that $\mL(z)=\Ps f$ for $f\in\ellP$ on $\oinftyZ$, 
then $z\equiv0$, cf. \cite[pg.~666]{amK89:I}. In other words, we assume the definiteness of system \Slaf{0}{f} for 
every $f\in\ellP$. In addition, this assumption implies that system \eqref{Sla} is definite as it follows from the 
choice $f\equiv0$,  by Theorem~\ref{T:3.6}(iii) and Corollary~\ref{C:3.6A}. This assumption establishes the density 
of $\dom\Tnod$ in $\tellP$, cf. \cite[Theorem~7.6]{amK89:I}. As noted in \cite{hB13}, a similar assumption is also 
needed for operators associated with system \eqref{E:dLHS} in \cite{yS06}.
\end{remark}

\begin{theorem}\label{T:dense.domain}
 If system \Slaf{0}{f} is definite on $\oinftyZ$, then $\dom\Tnod$ is dense in $\tellP$.
\end{theorem}
\begin{proof}
 Assume that $\dom\Tnod$ is not dense in $\tellP$. Then there exists  $\tf\in\dom(\Tnod)^\bot$ such that $\normP{f}\ne 0$.  
 Let $\tz\in\dom\Tnod$ be such that $\mL(z)=\Ps f$ and $\ty\in\dom\Tnod$ be such that $\mL(y)=\Ps g$ for some $g\in\ellP$.
 Then, by the extended Lagrange identity from Theorem~\ref{T:1},  we obtain
  \begin{equation*}
   \innerP{g}{z}-\innerP{y}{f}=\lim_{k\to\infty} y_{k+1}^*\,\J\,z_{k+1}-y_{0}^*\,\J\,z_{0}=0,
  \end{equation*}
 because $y,z\in\ellPn$, i.e., we have 
  \begin{equation}\label{E:thm.dens}
   \innerP{g}{z}=\innerP{y}{f}=0.
  \end{equation}
 Since $g\in\ellP$ was chosen arbitrarily and $z\in\ellPn$, we can take $g=z$ and the solution of $\mL(y)=\Ps z$ can be 
 obtained as $y_k=\Phk\,\J\sum_{j=0}^{k-1}\Ph_j^*\,\Ps_j\,z_j$. Then, \eqref{E:thm.dens} implies that
 $\innerP{z}{z}=0$, i.e., $\Psk\,z_k=0$ on $\oinftyZ$. Thus we have $\mL(z)=\Ps f$ and $z\in\tilde{0}$, which yields 
 $z\equiv0$ by the definiteness assumption for system \Slaf{0}{f}. Then $\Psk\,f_k=0$ on $\oinftyZ$, i.e., 
 $f\in\tilde{0}$, and the density of $\dom\Tnod$ in $\tellP$ is thus established.
\end{proof}

%%%%%%%%%%%%%%%%%%%%%%%%%%%%%%%%%%%%%%%%%%%%%%%%%%
%%%%%%%%%%%%%%%%%%%%%%%%%%%%%%%%%%%%%%%%%%%%%%%%%%
\subsection{The orthogonal decomposition of sequence spaces} 

Next, we introduce a linear map which will allow an orthogonal decomposition of $\tellPI{(\mI)}$, where $\mI\subset\oinftyZ$
denotes a discrete finite interval.  First, note that the sum $\sum_{k\in\mI}\Phk^*(\bla)\Psk u_k$ is independent of
$u\in\tz\in\tellPI{(\mI)}$, and let $ \K_{\la,\mI}$ denote the linear map defined by
 \begin{equation*}
 \K_{\la,\mI}:\tellPI{(\mI)}\to\Cbb^{2n},\qquad \K_{\la,\mI}(\tz)\coloneq \sum_{k\in\mI}\Phk^*(\bla)\Psk\zk,
 \end{equation*}
which we will abbreviate as $\K_{\mI}$ when $\la=0$.

\begin{lemma}\label{L:3.8}
 Let $\mI\subset\oinftyZ$ denote a  finite discrete subinterval, and let  $\la\in\Cbb$.  Then, $\ran\K_{\la,\mI}$ is 
 independent of $\la\in\Cbb$; in particular,
  \begin{equation}\label{E:3.8.1}
   \ran\K_{\la,\mI}=\{\xi\in\Cbb^{2n}\mid\ \Psk\,\Phk(\bla)\,\xi=0,\ \forall \ k\in\mI\}^\bot=\Ran\vp(\mI).
  \end{equation}
 Furthermore, $\tellPI{(\mI)}$ admits the following orthogonal sum decomposition:
  \begin{equation}\label{E:3.8.2}
   \tellPI{(\mI)}=\ker\K_{\la,\mI}\oplus\big\{ \tz\in\tellPI{(\mI)}\mid z=\Ph(\bla)\,\xi, \ \xi\in\Ran\vp(\mI)\big\}.
  \end{equation}
\end{lemma}
\begin{proof}
 For any $\xi\in\Cbb^{2n}$ and $\tu\in\tellPI{(\mI)}$,
  \begin{equation*}
   \inner{\xi}{\K_{\la,\mI}(\tu)}_{\Cbb^{2n}}=\sum_{k\in\mI}\xi^*\,\Phk^*(\bla)\,\Psk\,u_k
     =\innerPI{\mI}{\Ph(\bla)\,\xi}{u}.
  \end{equation*}
 Hence, $\K_{\la,\mI}^*:\Cbb^{2n}\to\tellPI{(\mI)}$, where  $\K_{\la,\mI}^*(\xi)=\tz$, and $z=\Ph(\bla)\,\xi$,  
 for $\xi\in\Cbb^{2n}$. In particular,
  \begin{equation*}
   \ran\K^*_{\la,\mI}=\{ \tz\in\tellPI{(\mI)}\mid \  z=\Ph(\bla)\,\xi, \   \xi\in\Cbb^{2n}   \},
  \end{equation*}
 and
  \begin{equation*}
   \ker\K_{\la,\mI}^*= \{\xi\in\Cbb^{2n}\mid \ \normPI{\mI}{\Ph(\bla)\,\xi}=0\}
   =\{\xi\in\Cbb^{2n}\mid \Psk\Phk(\bla)\xi =0, \forall \ k\in\mI\}.
  \end{equation*}
 Thus, the first of the equalities in \eqref{E:3.8.1} follows from the fact that
 $\ran\K_{\la,\mI}=\ker(\K_{\la,\mI}^*)^\bot$.

 Next, note that $\tellPI{(\mI)} = \ker\K_{\la,\mI} \oplus \ker(\K_{\la,\mI})^{\bot} = \ker\K_{\la,\mI}
 \oplus \ran\K^*_{\la,\mI}$. Let $\xi\in\Cbb^{2n}$, and $\xi=\eta+\zeta$, where $\eta\in\ran\K_{\la,\mI}$,
 and $\zeta\in\ran(\K_{\la,\mI})^\bot=\ker\K^*_{\la,\mI}$. Then, let $z=\Ph(\bla)\,\xi$, $u=\Ph(\bla)\,\eta$, 
 and $v=\Ph(\bla)\, \zeta$, and note that $\normPI{\mI}{v}=\normPI{\mI}{\Ph(\bla)\, \zeta}=0$. Thus, 
 $\tilde v=\tilde0\in\tellPI{(\mI)}$, and as a consequence, $\tz = \tu$, and hence
  \begin{equation*}
   \ran\K^*_{\la,\mI}=\{ \tz\in\tellPI{(\mI)}\mid \ z=\Ph(\bla)\,\eta, \ \eta\in\ran\K_{\la,\mI} \}.
  \end{equation*}
 Thus, to complete the demonstration of \eqref{E:3.8.1} and \eqref{E:3.8.2}, it remains to show that
 $\ran\K_{\la,\mI}=\Ran\vp(\mI)$.

 Let $\tz\in\tellPI{(\mI)}$. Then $\tz= \tv + \tu$, where $\tv\in\ker\K_{\la,\mI}$, and $u=\Ph(\bla)\eta$, where
 $\eta\in\ran\K_{\la,\mI}$. Then,
  \begin{equation*}
   \K_{\la,\mI}(\tz)=\K_{\la,\mI}(\tu)=\sum_{k\in\mI}\Phk^*(\bla)\,\Psk\,\Phk(\bla)\,\eta=\vp(\la,\mI)\,\eta,
  \end{equation*}
 and hence, $\ran\K_{\la,\mI}\subseteq\Ran\vp(\mI)$.  On the other hand, if $\xi\in\Cbb^{2n}$ and $z=\Ph(\bla)\xi$, then
  \begin{equation*}
   \vp(\la,\mI)\xi=\sum_{k\in\mI}\Phk^*(\bla)\,\Psk\,\Phk(\bla)\,\xi
     =\sum_{k\in\mI}\Phk^*(\bla)\,\Psk\,z_k=\K_{\la,\mI}(\tz),
  \end{equation*}
 thus showing that $\ran\K_{\la,\mI}=\Ran\vp(\mI)$.
\end{proof}

In analogy with $\K_{\la,\mI}$, we define the linear map
 \begin{equation*}
  \K_{\la}:\tellPo\to\Cbb^{2n},\qquad   \K_{\la}(\tz)\coloneq \sum_{k\in\Nbb_0}\Phk^*(\bla)\,\Psk\,\zk,
 \end{equation*}
which we abbreviate as $\K$ when $\la=0$. We note that the sum $\sum_{k\in\Nbb_0}\Ph^*_k(\bla)\Psk u_k$ is independent of 
$u\in\tz\in\tellPo$.

For the next two lemmas, recall from Lemma~\ref{L:5} (viz. \eqref{E:3.6}) that there is a discrete finite 
interval $\mI\subset\Nbb_0$ for which $\rank\vp(\mI)$ is maximal; that is, $\rank\vp(\mhI)=\rank\vp(\mI)$, for any 
discrete finite interval $\mhI$ such that $\mI\subseteq\mhI$. By \eqref{E:3.8.1}, it also follows that if $\mI$ is a 
discrete finite interval for which $\rank\vp(\mI)$ is maximal, then $\ran\K_{\la,\mhI}=\ran\K_{\la,\mI}$, when $\mhI$ is 
a discrete finite interval for which $\mI\subseteq\mhI$.

\begin{lemma}\label{L:3.9}
 Let $\mI\subset\Nbb_0$ be a discrete finite interval for which $\rank\vp(\mI)$ is maximal. Then, 
 $\ran\K_{\la}=\Ran\vp(\mI)$, for all $\la\in\Cbb$; in particular, $\ran\K_{\la}$ is independent of $\la$.
\end{lemma}
\begin{proof}
 Let $\tz\in\dom\K_\la=\tellPo$. Then, there exists an $N\in\Nbb_0$ such that $\Psk\zk=0$ for $k\ge N$; 
 hence, $z\in\ell^2_{\Psi}([0,N]_\Zbb)$, and $\tz\in\tellPI{([0,N]_\sZbb)}=\dom\K_{\la,[0,N]_\sZbb}$. In particular,
 $\K_{\la,[0,N]_\sZbb}(\tz)=\K_\la(\tz)$, and thus, $\ran\K_\la\subseteq\ran\K_{\la,[0,N]_\sZbb}$. Without loss of 
 generality, we may assume that $\mI\subseteq [0,N]_\sZbb$, and hence that 
 $\ran\K_{\la,\mI}=\ran\K_{\la,[0,N]_\sZbb} \supseteq \ran\K_\la$.

 Conversely, suppose that $\tz\in\tellPI{(\mI)}=\dom\K_{\la,\mI}$. Let $ u\in\ellPo$ be defined by
  \begin{equation*}
   u_k=\begin{cases}   
        z_k, & k\in\mI,\\  
        0, & k\not\in\mI.
       \end{cases}
  \end{equation*}
 Then, $\K_{\la,\mI}(\tz)= \sum_{k\in\mI}\Phk^*(\bla)\Psk z_k=\sum_{k\in\Nbb_0}\Phk^*(\bla)\Psk u_k=\K_\la(\tu)$. 
 Hence, $\ran\K_{\la,\mI}\subseteq\ran\K_\la$, thus $\ran\K_\la=\ran\K_{\la,\mI}$; which by Lemma~\ref{L:3.8} 
 yields $\ran\K_\la=\Ran\vp(\mI)$.
\end{proof}

\begin{lemma}\label{L:3.10}
 Let $\mI\subset\Nbb_0$ be a discrete finite interval for which $\rank\vp(\mI)$ is maximal. Then,
  \begin{equation*}
   \tellPo=\ker\K_\la \oplus \big\{ \tz\in\tellPo\mid z=\Ph(\bla)\,\xi, \ \xi\in\Ran\vp(\mI)\big\}
  \end{equation*}
 and $\codim(\ker\K_\la)=\rank\vp(\mI)$.
\end{lemma}
\begin{proof}
 In the complete analogy with the argument given in the proof of Lemma~\ref{L:3.8}, one can show that
  \begin{equation*}
   \tellPo= \ker\K_\la\oplus\ker(\K_\la)^\bot=\ker\K_\la\oplus\ran\K^*_\la,
  \end{equation*}
 where
  \begin{equation*}
   \ran\K^*_\la=\big\{ \tz\in\tellPo\mid z=\Ph(\bla)\,\xi, \ \xi\in\Ran\vp(\mI)\big\},
  \end{equation*}
 and that
  \begin{equation*}
   \ran(\K_\la)^\bot=\ker\K^*_\la=\{\xi\in\Cbb^{2n}\mid \ \normP{\Ph(\bla)\,\xi}=0\}
     =\{\xi\in\Cbb^{2n}\mid \Psk\Phk(\bla)\xi =0, \forall \ k\in\Nbb_0\}.
  \end{equation*}

 Let $\xi_j,\ j=1,\dots, m$, denote a basis for $\ran\vp(\mI)=\ran\K_\la$, let $u^{[j]}\coloneq\Phi(\bla)\xi_j$, and thus
 $\tu^{[j]}\in\ran\K^*_\la$, for $ j=1,\dots,m$. Suppose that $\sum_{j=1}^{m}\alpha_j\tu^{[j]}=\tilde 0\in\tellPo$. Then, 
 $\Psk\Phk(\bla)(\sum_{j=1}^m\alpha_j\xi_j )=0$, for all $k\in\Nbb_0$. As a consequence,
  \begin{equation*}
   \sum_{j=1}^m \alpha_j\xi_j \in \ran\K_\la \cap \ran(\K_\la)^\bot=\{0\},
  \end{equation*}
 and hence, $\alpha_j=0,\ j=1,\dots, m$; thus, $\tu^{[j]}\in\ran\K^*_\la,\ j=1,\dots,m$, are independent 
 in $\tellPo$. Hence, $\codim(\ker\K_\la)=\dim(\ran(\K^*_\la))=\rank\vp(\mI)$.
\end{proof}

%%%%%%%%%%%%%%%%%%%%%%%%%%%%%%%%%%%%%%%%%%%%%%%%%%%%%%%%
%%%%%%%%%%%%%%%%%%%%%%%%%%%%%%%%%%%%%%%%%%%%%%%%%%%%%%%%

\subsection{The minimal relation  and its deficiency indices}

Before we define the minimal linear relation $\Tmin$ and establish the fundamental relation between $\Tmax$ and 
$\Tmin$, we prove two auxiliary lemmas.

\begin{lemma}\label{L:4.1}
 We have $\ker\K_{\la}\subseteq\ran(\Tnod-\la I )$ for every $\la\in\Cbb$.
\end{lemma}
\begin{proof}
 Let $\la\in\Cbb$ and $\tf\in\ker\K_\la\subseteq\tellPo$. Then, there is $N\in\Nbb_0$ such that
 $\Psk g_k=0$ for all $k\ge N$ and all $g\in\tf$. Moreover, $\K_\la(\tf)=\sum_{k\in\Nbb_0}\Phk^*(\bla)\,\Psk\,g_k=0$  for
 all $g\in\tf$.

 For $g\in\tf$, we note that 
  \begin{equation}\label{E:4.4A}
   z_k\coloneq-\Phk(\la)\J\sum_{j=k}^\infty\Phi_j^*(\bla)\Psi_jg_j=
    \begin{cases}
      -\Phk(\la)\J\sum_{j=k}^{N-1}\Phi_j^*(\bla)\Psi_jg_j, & k\in[0,N-1]_\Zbb,\\
      0,& k\in[N,\infty)_\Zbb,
    \end{cases}
  \end{equation}
 defines a sequence $z\in\ellP$ such that $\mL(z)=\la\Psi z+\Psi g$ (cf. \cite[Theorem~3.17]{sE05})
 for which $z_0=-\Phi_0(\la)\J\K_\la(\tf)=0$. Thus, as defined, $z\in\ellPn$ and satisfies
 $\mL(z)-\la\Psi z= \Psi g$; thereby showing that $\tf\in\ran(T_0-\la I)$; cf. \eqref{e5.3}.
\end{proof}

\begin{lemma}\label{L:4.2}
 We have $\innerP{\tf}{\ty}=\innerP{\tz}{\tg}$\ for  every $\{\tz,\tf\}\in\Tnod$ and $\{\ty,\tg\}\in\Tmax$.
\end{lemma}
\begin{proof}
 Let $\{\tz,\tf\}\in\Tnod$ and $\{\ty,\tg\}\in\Tmax$. Let $z\in\ellPn$ such that $\mL(z)=\Psi f$,
 and let $y\in\ellP$ such that $\mL(y)=\Psi g$. Then, by the extended Lagrange identity in \eqref{E:8} 
 with $\la=\nu=0$, we see that
  \begin{equation*}
   \innerPI{[0,k]_\sZbb}{z}{g}= -z^*_j\J y_j\big|_0^{k+1} + \sum_{j=0}^k f^*_j\Psi_j y_j
     =-z^*_j\J y_j\big|_0^{k+1} + \innerPI{[0,k]_\sZbb}{f}{y}.
  \end{equation*}
 However, $z\in\ellPn$ implies that $z_0=0$ and that $z_k=0$ for sufficiently large $k\in\Nbb_0$. As
 a consequence we see that $\innerP{z}{g}=\innerP{f}{y}$, and hence, $\innerP{u}{g}=\innerP{f}{w}$
 for all $u\in\tz$, and $w\in\ty$.
\end{proof}

By Lemma~\ref{L:4.2} and \eqref{E:2.2.1A}, one obtains
 \begin{equation}\label{E:4.6}
  \Tnod\subseteq\Tmax\subseteq\Tnod^*,
 \end{equation}
from which we conclude that $\Tnod$ is symmetric in $\tellP^2$, and hence closable. We then define
the minimal operator, $\Tmin$, by
 \begin{equation*}
  \Tmin\coloneq \overline{\Tnod}.
 \end{equation*}
Another approach to defining a minimal operator is to let $\Tmin\coloneq\Tmax^*$; cf.
\cite[Definition~2.3]{mL.mmM03} and \cite[Identity~(4.2)]{jB.sH.hsvdS.rW11}. In the examples cited, this is
equivalent to our choice of definition for $\Tmin$. The next theorem demonstrates this for the
operators at hand. We note also that an alternative demonstration of the next result can be
patterned after that presented in \cite[pp.~1354--1355]{jB.sH.hsvdS.rW11} using the results in
Section~\ref{S:nonhom} and \cite[Proposition~A.2]{jB.sH.hsvdS.rW11}.

\begin{theorem}\label{T:4.3}
 For $\Tmax$ and $\Tnod$ as defined in \eqref{E:4.1} and \eqref{E:4.2}, respectively,
  \begin{equation}\label{E:4.7}
   \Tnod^*=\Tmin^*=\Tmax.
  \end{equation}
\end{theorem}
\begin{proof}
 By \eqref{E:4.6}, note that $\Tmax\subseteq\Tnod^*=\overline{\Tnod}^*=\Tmin^*$. Thus, it
 remains to show that $\Tnod^*\subseteq\Tmax$; or equivalently, given $\{\tz,\tf\}\in\Tnod^*$, that 
 there is a $v\in\tz\in\tellP$ such that $\mL(v)=\Psi f$.

 Now, let $\{\tu,\tg\}\in\Tnod$ and note that there exists $y\in\tu$ such that $y\in\ellPn$ and $\mL(y)_k=\Psk\gk$ 
 for $k\in\oinftyZ$. Next let $w\in\Cbb(\Nbb_0)^{2n}$ be a sequence such that $\mL(w)=\Psi f$, where 
 $\tf\in\Tnod^*(\tz)$. Then, by the Lagrange identity \eqref{E:8}, we see that 
 $\innerPI{[0,k]_\sZbb}{w}{g}=w^*_j\J y_j\big|_0^{k+1} + \innerPI{[0,k]_\sZbb}{f}{y}$. Given that $y_0=0$ and the $y_k=0$ 
 for sufficiently large $k\in\Nbb_0$, we see that $\innerP{w}{g}=\innerP{f}{y}$. From the definition of 
 $\Tnod^*$, note that $\innerP{\tz}{\tg}=\innerP{\tf}{\tu}$.  As a consequence, $\innerP{z-w}{g}=0$ for any 
 $z\in\tz\in\dom \Tnod^*$, and any $g\in\tg\in\ran T_0$. Recall from Lemma~\ref{L:4.1}, with $\la=0$, that
 $\ker\K\subseteq\ran T_0\subseteq\tellPo$. Then, for every $h\in\tilde{h}\in\ker\K$ we see that $\innerP{z-w}{h}=0$, 
 and that $\sum_{k\in\Nbb_0}\Phk^*\Psk\,h_k=0$ (i.e., $\Phk=\Phk(0)$); and as a~consequence, for any $\xi\in\Cbb^{2n}$, 
 that
  \begin{equation}\label{E:4.10}
   \sum_{k=0}^\infty(z_k-w_k-\Phk\,\xi)^*\,\Psk\,h_k=0.
  \end{equation}

 Next, let $\mI\subset\Nbb_0$ be a finite discrete interval such that $\rank\vp(\mI)=m\le 2n$ is
 maximal as noted in Lemma~\ref{L:5}. Then, by Lemma~\ref{L:3.10}, we see that $\tellPo=\ker\K\oplus\ker\K^\bot$ 
 and that there is a~basis, $\ty^{[1]},\dots,\ty^{[m]}\in\tellPo$, for $\ker\K^\bot$, in which 
 $y^{[j]}=\Phi\xi^{[j]}$, and where $\xi^{[j]}, j=1,\dots,m$, forms a~basis for $\ran\vp(\mI)$.
 Then, there is discrete finite interval $\mhI=[0,N]_\Zbb$, where $\mI\subseteq\mhI$, for
 which
  \begin{equation}\label{E:4.11}
   \Psk\,y_k^{[j]}=0,\quad k\ge N,\quad j=1,\dots,m.
  \end{equation}
 Now, suppose that $\tg\in\ker\K_{\mhI}$, then, for $j=1,\dots, m$,
  \begin{equation*}
   \innerPI{\mhI}{y^{[j]}}{g}=\sum_{k\in\mhI}{\xi^{[j]}}^*\Phk^*\Psk g_k={\xi^{[j]}}^*\K_{\mhI}(\tg)=0;
  \end{equation*}
 and thus, $\ty^{[j]}\in\ker(\K_{\mhI})^\bot$, for $j=1,\dots, m$.

 Let $x\in\ellPo$ be defined by
  \begin{equation*}
   x_k=\begin{cases} z_k-w_k, & k\in\mhI,\\ 0, & k\in\Nbb_0\setminus\mhI,\end{cases}
  \end{equation*}
 where $\{\tz,\tf\}\in\Tnod^*$ and $w\in\Cbb(\Nbb_0)^{2n}$ with $\mL(w)=\Psi f$. Now, 
 $\tilde x\in\tellPI{(\mI)}$ and, by Lemma~\ref{L:3.8}, there is $\tv\in\ker(\K_{\mhI})^\bot$, 
 where $v=\Phi\eta$ for some $\eta\in\Ran\vp(\mhI)$, such that $\tilde x-\tv\in\ker\K_{\mhI}$. 
 As a~consequence $\innerPI{\mhI}{\tilde x - \tv}{\ty^{[j]}}=0$ for all $j=1, \dots, m$. Then, by
 \eqref{E:4.11}, we obtain, for $j=1,\dots,m$, that
  \begin{equation}\label{e5.11}
   \sum_{k\in\Nbb_0}(z_k-w_k-\Phk\eta)^*\Psk y^{[j]}_k=\sum_{k\in\mhI}(z_k-w_k-\Phk\eta)^*\Psk y^{[j]}_k
     =\innerPI{\mhI}{\tilde x - \tv}{\ty^{[j]}}=0.
  \end{equation}
 Together, \eqref{E:4.10} and \eqref{e5.11} show that $\innerP{z-w-\Phi\eta}{h}=0$ for all $h\in\ellPo$; 
 in particular, for
  \begin{equation*}
   h_k=\begin{cases} z_{k_1}-w_{k_1}-\Phi_{k_1}\eta, & k=k_1,\\ 0, & k\in\Nbb_0\setminus\{k_1\}.\end{cases}
  \end{equation*}
 Thus, for any $k_1\in\Nbb_0$, we see that $\Psi_{k_1}(z_{k_1}-w_{k_1}-\Phi_{k_1}\eta)=0$, 
 and hence $\normP{z-(w+\Phi\eta)}=0$. In summary, we get that $w+\Phi\eta\in\tz$ where
 $\{\tz,\tf\}\in\Tnod^*$, and further that $\mL(w+\Phi\eta)=\mL(w)=\Ps f$; hence, showing that
 $\Tnod^*\subseteq\Tmax$.
\end{proof}

Recalling the definition of deficiency subspaces and indices from Section~\ref{S:2.2} 
(viz. \eqref{e2.17}, \eqref{e2.18}), let 
 \begin{gather*}
  \Mla\coloneq \{z\in\ellP\mid \mL(z)_k=\la\Psk\,\zk,\ k\in\oinftyZ\},
   \quad \dla\coloneq \dim \Mla,\\ 
  \tMla\coloneq \big\{\tz\in\tellP\mid\{\tz,\la\tz\}\in\Tmax=\Tmin^*\big\},\quad \tdla\coloneq \dim \tMla.
 \end{gather*}
In other words, the subspace $\Mla\subseteq\ellP$ consists of all square summable solutions of
\eqref{Sla}, while $\dla$ denotes the number of linearly independent square summable solutions of 
\eqref{Sla}. Then, $\tMla=\tMla(\Tmin)$ represents the deficiency space, and $\tdla=\tdla(\Tmin)$ the 
deficiency index, for $\Tmin$ and $\la\in\Cbb$. It was shown in Theorem~\ref{T:4.3} that $\Tmin$ 
is a closed, symmetric linear relation, and hence $\tdla$ is constant on each of the open upper and lower 
half-planes of $\Cbb$; cf. \cite[Theorem~2.13]{yS12}. Thus we let $\td_+\coloneq\td_{i}=\tdla$, for 
$\la\in\Cbb_+$, and $\td_-\coloneq\td_{-i}=\tdla$, for $\la\in\Cbb_-$. Then, with the aid of
\cite[Proposition~A.1]{jB.sH.hsvdS.rW11} we obtain the following von~Neumann decomposition of the maximal 
linear relation $\Tmax$ in terms of the minimal linear relation $\Tmin$ and the defect subspaces $\tMla$ 
and $\tMbla$, i.e.,
 \begin{equation*}
  \Tmax=\Tmin\dotplus\widetilde{\mathcal{M}}_\la\dotplus\widetilde{\mathcal{M}}_{\bla},
 \end{equation*}
where $\widetilde{\mathcal{M}}_\la\coloneq\big\{\{\tz,\la\tz\}\mid\{\tz,\la\tz\}\in\Tmax\big\}$ and the direct sum $\dotplus$ is 
orthogonal if $\la=\pm i$.

Assuming that $\rank\vp(\mI)$ is maximal, as described in Lemma~\ref{L:5}, we next show a relationship
between $\dla$ and $\tdla$.

\begin{theorem}\label{T:4.5}
 Let $\mI\subset\Nbb_0$ be a discrete finite interval such that $\rank\vp(\mI)$ is maximal. Then, for any $\la\in\Cbb$,
  \begin{equation}\label{E:4.20}
   \dla=\tdla+2n-\rank\vp(\mI).
  \end{equation}
\end{theorem}
\begin{proof}
 Let $\la\in\Cbb$, and $\tz\in\tMla$. Then there exists $z\in\tz$ which solves system \eqref{Sla}, i.e.,
 $z\in\Mla$. With $\pi_1$ denoting the restriction of the quotient space map of Section~\ref{SS:2.2} given by  
 $\pi_1\coloneq \pi\vert_{\Mla}:\Mla\to\tMla$, we note that
  \begin{equation*}
   \ker\pi_1=\{z\in\ellP\mid\ \mL(z)_k=\la\Psk\,\zk\text{ and }\Psk\,\zk=0 \text{ for } k\in\oinftyZ\}.
  \end{equation*}
 Then $z\in\ker\pi_1$ if and only if $z=\Phi\xi$ for some $\xi\in\bigcap_{n\in\Nbb_0}\ker\vp([0,n]_\Zbb)=\ker\vp(\mI)$.  
 Hence $\dim (\ker\pi_1)=\dim(\Ker\vp(\mI))=2n-\rank\vp(\mI)$. Since the map $\pi_1$ is surjective,
  \begin{equation*}
   \dim\Mla=\dim\tMla+\dim(\Ker\pi_1),
  \end{equation*}
 which implies \eqref{E:4.20}.
\end{proof}

The following properties for $\dla$ and $\tdla$ are a direct consequence of \eqref{E:4.20}.

\begin{corollary}\label{C:4.6}
 Let $\la\in\Cbb$. Then,
  \begin{itemize}
    \item[\rm(i)] $\tdla=\dla$ if and only if system~\eqref{Sla} is definite on $\oinftyZ$;
    \item[\rm(ii)] $\dla-\tdla$ is nonnegative and constant for $\la\in\Cbb$;
    \item[\rm(iii)] $\dla$ is constant in $\Cbb_+\coloneq \{\la\in\Cbb\mid\ \Im(\la)>0\}$ and in
                    $\Cbb_-\coloneq \{\la\in\Cbb\mid\ \Im(\la)<0\}$.
  \end{itemize}
\end{corollary}
\begin{proof}
 The first statement follows directly from \eqref{E:4.20} and Theorem~\ref{T:3.6}, while
 the second statement is a consequence of \eqref{E:4.20} and Lemma~\ref{L:3.4}.  Since  the linear
 relation $\Tnod$ is symmetric (cf. \eqref{E:4.6}), the number $\tdla$ is constant in $\Cbb_+$ and 
 $\Cbb_-$ by \cite[Theorem~2.13]{yS12}. Hence, the value of $\dla$ is also constant in $\Cbb_+$ and $\Cbb_-$  
 by Theorem~\ref{T:4.5}, and by the independence of $\rank\vp(\mI)$ on the value of  $\la$ established
 in Lemma~\ref{L:3.4}.
\end{proof}

\begin{remark}
 The last statement of Corollary~\ref{C:4.6} extends the enumeration and analysis of linearly independent 
 square summable solutions  for system \eqref{Sla} found in \cite{rSH.pZ:W-TGLP} 
 (viz. \cite[Theorem~4.9]{rSH.pZ:W-TGLP}) by showing that $r(\la)$ defined in \cite[Identity~(4.1)]{rSH.pZ:W-TGLP} 
 is constant in $\Cbb_+$ and $\Cbb_-$.
\end{remark}

\begin{theorem}\label{T:4.7}
 For $\la\in\Cbb$ we have
  \begin{equation}\label{e5.13}
   \ker(\Tmin-\la I)=\{\tilde 0\}.
  \end{equation}
 Moreover, for every $\la\in\Rbb$ we have $\td_{\la}\leq\td_{\pm}$.
\end{theorem}
\begin{proof}
 First, note that $\tellPo$ is dense in $\tellP$. Then, for $\la\in\Cbb$ and  $\tf\in\tellPo$, and with 
 $z\in\ellPo$ as defined in \eqref{E:4.4A}, we note, by the construction given in Lemma~\ref{L:4.1}, that
 $\{\tz,\tf\}\in\Tmax-\la I$. This shows  that $\tellPo\subseteq\ran(\Tmax -\la I)$ and hence that 
 $\ran(\Tmax - \la I)$ is dense in $\tellP$ for any $\la\in\Cbb$. As a consequence of \eqref{E:4.7}, we get
 that $\ker(\Tmin-\la I)=\ran(\Tmax-\bla I)^\bot=\{\tilde 0\}$.

 The remaining assertion follows from \cite[Lemma~2.25]{mL.mmM03} and \eqref{e5.13}, as noted in the comments 
 associated with \eqref{E:2.2.3} at the end of Section~\ref{S:2.2}.
\end{proof}

By Corollary~\ref{C:3.6A}, Corollary~\ref{C:4.6}(i), and Theorem~\ref{T:4.7} we obtain the following statement.

\begin{corollary}
 When system~\eqref{Sla} is definite on $\oinftyZ$, then $d_{\la}\leq d_{\pm}$ for any $\la\in\Rbb$.
\end{corollary}

%%%%%%%%%%%%%%%%%%%%%%%%%%%%%%%%%%%%%%%%%%% SECTION %%%%%%%%%%%%%%%%%%%%%%%%%%%%%%%%%%%%%%%%%%%%%%%%%%%%%%%%%%%%%%%%%%

\section*{Acknowledgements}
This work was supported by the Program of ``Employment of Newly Graduated Doctors of Science for Scientific Excellence'' 
(grant number CZ.1.07/2.3.00/30.0009) co-financed from European Social Fund and the state budget of the Czech Republic. The 
second author would like to express his thanks to the Department of Mathematics and Statistics (Missouri University of Science 
and Technology) for hosting his visit. The authors are also indebted to the anonymous referees for detailed reading of the 
manuscript and constructive comments in their reports which helped to improve the presentation of the results.

%%%%%%%%%%%%%%%%%%%%%%%%%%%%%%%%%%%%%%%%%%% BIBLIOGRAPHY %%%%%%%%%%%%%%%%%%%%%%%%%%%%%%%%%%%%%%%%%%%%%%%%%%%%%%%%%%%%%
% \bibliographystyle{my_bibliography_paper_style}
% \bibliography{bibliotheca}

\end{document}